
\magnification=\magstep1

\input amssym.def
\input amssym.tex
\font\bigbf=cmbx12
%
%
\font\teneusm=eusm10
\font\seveneusm=eusm7
\font\fiveeusm=eusm5
\newfam\eusmfam
\textfont\eusmfam=\teneusm
\scriptfont\eusmfam=\seveneusm
\scriptscriptfont\eusmfam=\fiveeusm

\font\tenmib=cmmib10
\font\sevenmib=cmmib7
\font\fivemib=cmmib5
\newfam\mibfam
\textfont\mibfam=\tenmib
\scriptfont\mibfam=\sevenmib
\scriptscriptfont\mibfam=\fivemib

\font\tenss=cmss10
\font\sevenss=cmss8 scaled 833
\font\fivess=cmr5
\newfam\ssfam
\textfont\ssfam=\tenss
\scriptfont\ssfam=\sevenss
\scriptscriptfont\ssfam=\fivess
\def\ss{\fam\ssfam}
\thinmuskip = 2mu
\medmuskip = 2.5mu plus 1.5mu minus 2.1mu  
\thickmuskip = 4mu plus 6mu
\def\loosegraf#1\par{{%
\baselineskip=13.4pt plus 1pt \lineskiplimit=1pt \lineskip=1.3 pt
#1\par}}

\def\cc{{\Bbb C}}
\def\pp{{\Bbb P}}
\def\eps{\varepsilon}

\def\rr{{\Bbb R}}
\def\ss{{\Bbb S}}

\def\bb{{\Bbb B}}

\centerline{\bigbf One example in concern with extension and separate
analyticity }
\centerline{\bigbf properties of meromorphic mappings}
\bigskip\rm
\centerline{\rm S.~Ivashkovich\footnote{}{AMS subject classification:
32 D 15, Key words: meromorphic map, Hartogs extension theorem, Lelong number,
 separate meromorphi\-city, Rothstein-type extension theorem.
}}
\smallskip
\centerline{March 1994}

\bigskip
\centerline{\bf Abstract}
\def\longpoints{\leaders\hbox to 0.5em{\hss.\hss}\hfill \hskip0pt}
\medskip
\centerline{
\vbox{\hsize=4truein
\baselineskip=10pt \sevenrm
We construct a compact complex manifold of dimension three, such  that
every meromorphic map from two-dimensional domain into this manifold
extends meromorphically onto the envelope of holomorphy of this domain,
but there is a meromorphic map of punctured three-dimensional ball into our
manifold, which doesn't extend to origin. A Theorem describing the
obstructions occuring here is given.
}}

\bigskip
\noindent
\bf 0. Introduction.
\medskip
\rm The purpose of this paper is to present an example of compact complex
three-fold having "strange" behavior with respect to a meromorphic
mappings into it and give certain positive results in entitled direction.

All complex spaces, which we consider in this paper are supposed to be reduced
and normal.

\rm Let  $D$ and  $X$ be complex spaces. Recall that a meromorphic mapping
$f$ from $D$
to $X$ is defined by its graph $\Gamma_f$, which is an analytic subset of the
product $D\times X$, satisfying the following conditions:

\rm (i) $\rm{\Gamma_f}$ is locally  irreducible analytic subset of $D\times X$.

(ii) The restriction $\pi \mid _{\Gamma_f} : \Gamma_f \longrightarrow D$ of
the natural projection $\pi : D\times X \longrightarrow D$ to $\Gamma_f$
is proper, surjective and generically one to one.

\rm This notion of meromorphicity is due to Remmert, see [Re], and is based on
the
observation that meromorphic functions on $D$ are precisely the meromorphic
mappings, in the sense just defined  into $X=\cc\pp^1$.
Denote by

$$
H^n(r) = \left\{ (z_1,...z_n) \in C^n:\Vert z'\Vert < r, \vert z_n\vert < 1
 \hbox{ or } \Vert z'\Vert < 1, 1 - r <\vert z_n\vert <1 \right\} \eqno(0.1)
$$
\medskip
\noindent
a $n$-dimensional Hartogs figure. Here $z'= (z_1,...,z_{n-1})$, $0<r<1$
and
$\Vert  \cdot  \Vert $
stands for the polydisk norm in $\cc^n$.

\smallskip\noindent\bf
Definition. \it Recall, that the complex space $X$ possesses
a meromorphic (holomorphic)  extension property in dimension $n$ if  any
meromorphic (holomorphic) mapping $f:H^n(r)\to X$ extends meromorphically
(holomorphically) onto the unit polydisk $\Delta^n$.

\smallskip\rm
When one studies an extension properties of holomorphis mappings one finds
the following statements, proved correspondingly in [Sh-1] and [Iv-1],
usefull:
\smallskip\noindent\it
1. If a complex space $X$ posseds a holomorphic extension property in dimension
$n$ then for every domain $D$ over a Stein manifold $\Omega $ of dimension
$n$ any holomorphic mapping $f:D\to X$ extends to a holomorphic mapping
$\hat f:\hat D\to X$ from the envelope of holomorphy $\hat D$ of $D$ into
$X$.

\noindent
2. If a complex space $X$ posseds a holomopric extension property in
dimension $2$ then $X$ posseds a hol.ext.prop. in all dimensions $n\ge 2$.
\smallskip
\rm The same proof as was given by B. Shiffman in [Sh-1] shows that the first
 statement remains valid
for the meromrophic mappings. While, and this is our starting point here, the
second fails to be true!

 Namely in \S 1  we shall construct the following
\smallskip
\noindent
\bf{Example }. \it There exists a compact complex three-fold $X$ such that:

(a) For every domain $D$ in $\cc^2$ every meromorphic mapping $f :
D\longrightarrow X$ extends to a meromorphic mapping $\hat f : \hat D
\longrightarrow X$. Here $\hat D$ stands for the envelope of holomorphy of
$D$.

(b) But there exists a meromorphic mapping $F : \bb^3 \setminus \{ 0\}
\longrightarrow X$ from punctured three-ball into $X$ which does not extend
to the origin.

\smallskip\rm

\rm This example shows the essential difference between extension
properties of holomorphic and meromorphic mappings. In fact this phenomena
backwards several difficulties in studying the meromorphic maps. Thats why
we devote the "positive" part of this paper to the explainations what are
the obstructions for passing from extension of mappings of two-dimensional
domains to three-dimensional ones. We shall prove in this direction the
following
\smallskip
\noindent
\bf Theorem. \it  Let $f: H^3(r) \longrightarrow X $ be a meromorphic
mapping of $3$-dimensional Hartogs figure into a  complex space $X$,
which possesses a meromorphic extension property in dimension two. Then:
\medskip
(i) $f$ extends to a meromorphic mapping of $\Delta ^3 \setminus S $,
where $S$ has a form

$S = S_1\times S_2\times S_3$, and $S_j$ are closed
subsets in $\Delta $ of harmonic measure zero.

\medskip
(ii) Let $s_0\in S$ be a point. Suppose that  the Lelong numbers
$\Theta (f^{\ast }w_h,s_0)$

and $\Theta ((f^{\ast }w_h)^2,s_0)$ are finite.
Then there exists a regular modification

$\pi : \hat \Delta^3
\longrightarrow \Delta^3$ , such that  $f\circ \pi$ holomorphically extends
 onto some

neighborhood of $\pi^{-1}(s_0)$.

\smallskip
\noindent
\bf Remark. \rm  In this theorem we suppose that the space $X$ is equipped
with
some Hermitian metric $h$, and $w_h$ its associated (1,1)-form. Currents
$T_k = (f^{\ast }w)^k$  are defined on $\Delta ^3 \setminus S $. We define
the Lelong numbers of (k,k)-current $R$ in $s_0$ as
\medskip
$$
\Theta (R,s_0) = \lim_{\varepsilon \searrow 0}\sup 1/\varepsilon
^{2(n-k)}\int_{B_{s_0}(\eps )\setminus S} R \wedge  (dd^c\Vert z \Vert ^2)^
{2(n-k)}.
\eqno(0.2)
$$
\medskip
\noindent
$\Theta (R,s_0) $ can take also infinite values. Remark only that finiteness
or infiniteness of those numbers in $s_0$ doesn't depend on the particular
choice of metric $h$ on $X$.

We want to point out here the following. If one wishes to use the Bishops
extension theorem for analytic sets in order to extend the graph $\Gamma _f$
from $(\Delta^3\setminus S)\times X$ onto $\Delta^3\times X$, one needs
to
estimate the volume of $\Gamma _f$ in the neighbourhood of $S$, say in the
neighbourhood of some point $s_0\in S$. For this one needs to estimate besides

$$
\int_{B_{s_0}(\varepsilon )}f^{\ast }w\wedge (dd^c\Vert z\Vert ^2)^2
\hbox { and }
 \int_{B_{s_0}(\varepsilon )}(f^{\ast }w)^2\wedge dd^c\Vert z\Vert^2
 $$
\noindent
also the integral $\int_{B_{s_0}(\varepsilon )}(f^{\ast }w)^3$.
The theorem
stated above makes possible for us to get rid of the last integral.

Of course the finitness of Lelong numbers of the current $f^*w_h$ and its
powers are nessessary for the extension of $f$. Finally, holomorphic extension
of $f\circ \pi $ to the neighborhood of $\pi^{-1}(s_0)$ implies a meromorphic
extension of $f$ to the neighborhood of $s_0$.

Methods which are developped here we apply
to the  questions of separate analyticity of meromrophic mappings with
values in general complex spaces, see 2.5 for details.

I would also like to give my thanks to the referee for numerous remarks
and corrections.

\bigskip\noindent\bf
Table of content.\medskip \sl
\noindent
\line{0. Introduction. \longpoints pp. 1-3}

\smallskip\noindent
\line{1. Counterexample. \longpoints pp. 3 - 9}

1.1. Construction of an Example . 1.2. Construction of the

universal cover and nonextendable map. 1.3. Two-dimensional

extension property of $X$.

\smallskip\noindent
\line{2.  Meromorphic polydisks. \longpoints pp. 9 - 20}

2.1. Generalities on pluripotential theory. 2.2. Passing to

higher dimensions: holomorphic case. 2.3. Sequences

of analytic sets of bounded volume. 2.4. Families of

meromorphic  polydisks. 2.5. Application: Rothstein-type theorem

and separate meromorphicity.
\smallskip\noindent
\line{3. Estimates of Lelong numbers from below. \longpoints pp. 20 - 29}

3.1. Generalities on blowings-up. 3.2. Estimates. 3.3. Proof of

the Theorem.

\smallskip\noindent
\line{4. References. \longpoints pp. 29 - 31}

\bigskip\noindent\bf
1. Counterexample.
\smallskip\rm

We shall start with the construction of the example announced in the
Introduction.

\smallskip
\noindent
{\sl 1.1. Construction of an example.}
\smallskip
\rm Take $\cc^3$ with coordinates $z_1, z_2, z_3$. Denote by $l_0$ the line
$\{ z_3 = z_2 = 0 \}$. $l_0$ has a natural coordinate $z_1$. By $X_0$ denote
the ball of radii $1/2$ centered at zero in $\cc^3$.

Consider now a blowing up $\hat \cc^3$ of $C^3$ along $l_0$. The origin in
$\cc^3$
 we denote by $0_0$. Let $\pi_0 :\hat \cc^3 \longrightarrow \cc^3$ be our
modification, $E_0 = \pi_0^{-1}(l_0)$ -- the exceptional divisor. By $0_1$ we
denote
the point of intersection of the strict transform  $\pi_0^{-1}(\lbrack z_3
\rbrack )$ of $z_3$ - axis with $E_0$. In the affine neighborhood of $0_1$ we
introduce the standard coordinate system $u_1,u_2,u_3$ such that modification
$\pi_0 : \hat \cc^3 \longrightarrow \cc^3$ in those coordinates is given by

$$z_1 = u_1$$
$$z_2 = u_2u_3           \eqno(1.1.1)$$
$$z_3 = u_3$$
\smallskip
The preimage of $0_0$ under the modification $\pi_0$ we denote by $l_1^{'}$.
Put $\pi_0^{-1}(X_0) = X_1$. Now we blow up $X_1$ along $l_1^{'}$. Denote by
$\pi_1 : X_2 \longrightarrow X_1$ this modification. Let $E_1$ be the
exceptional  divisor of $\pi_1$ and $l_1$ let be the intersection of $E_1$
and $E_0$ (more precisely with the strict transform of $E_0$ under $\pi _1 $)
. The proper preimage of $0_1$ under $\pi_1$ we denote by $l_2^{'}$. The point
 of intersection of $l_2^{'}$ with $E_0$  denote by $0_2$. In the affine
neighborhood of $0_2$ we introduce coordinate  system $v_1,v_2,v_3$ (standard
one for the modification $\pi_1$) in which our modification $\pi_1$ is given by

$$u_1 = v_1$$
$$u_2 = v_2             \eqno(1.1.2)$$
$$u_3 = v_1v_3$$
\smallskip
Next, consider the blowing up of $X_2$ along $l_2^{'}$   $\pi_2 : X_3
\longrightarrow X_2 $. By $E_2$ denote the exceptional divisor of $\pi_2$.
Put $l_2 = E_2\cap E_1, l_3^{'} = E_2\cap E_0$. $0_3$ let be the point of
intersection $E_0\cap E_1\cap E_2 $. Here again by $E_0,E_1$ we understand
the appropriate strict transforms. We introduce a coordinates $w_1,w_2,w_3$
in affine neighborhood of $0_3$ in such a way that $\pi_2 : X_3
\longrightarrow X_2$ is given in those coordinates by

$$v_1 = w_1w_2$$
$$v_2 = w_2    \eqno(1.1.3)$$
$$v_3 = w_3$$

Take now the composition $\pi = \pi_0\circ \pi_1\circ \pi_2: X_3
\longrightarrow X_0 $,which in coordinates can be written as

$$z_1 = w_1w_2$$
$$z_2 = w_1w_2^2w_3  \eqno(1.1.4)$$
$$z_3 = w_1w_2w_3$$

Take now two copies of $X_0$, denote them by $X_0^{(1)}$ and $X_0^{(2)}$. On
$X_0^{(i)}$ let us fix initial coordinates $z_1^{(i)},z_2^{(i)},z_3^{(i)}$.
On the affine neighborhood of $0_3^{(i)}\in X_3^{(i)}$ we fix coordinates
$w_1^{(i)},w_2^{(i)},w_3^{(i)}$ so as where introduced by $(1.1.1), (1.1.2),
(1.1.3)$.
Consider now a holomorphic mapping $\phi : X_0^{(1)} \longrightarrow
X_3^{(2)}$, which sends zero $0_0^{(1)}$ from $X_0^{(1)}$ to $0_3^{(2)}$ and
in coordinates $z^{(1)}$ and $w^{(2)}$ is defined by

$$w_1^{(2)} = z_1^{(1)}$$
$$w_2^{(2)} = z_2^{(1)}   \eqno(1.1.5)          $$
$$w_3^{(2)} = z_3^{(1)}$$

\noindent
i.e. $\phi $ is identity in $z^{(1)}, w^{(2)}$.

Now blow up $l_0^{(1)}$ in $X_0^{(1)}$ and $l_3^{(2)^{'}}$ in $X_3^{(2)}$.
Denote the manifolds obtained by $\hat X_0^{(1)}$ and $\hat X_3^{(2)}$.
$\hat X_0^{(1)} = X_1^{(1)}$ in our notations. Proper preimage of $l_3^{
(2)^{'}}$ denote by $E_3$.  Lift up $\phi $ to a biholomorphism $\hat
\phi : \hat X_0^{(1)} \longrightarrow \hat \phi (\hat X_0^{(1)}) \subset
\hat X_3^{(2)}$. The fact that $\phi $ lifts to a
biholomorphism, not just a bimeromorphism, of blowings up follows from the
observation that $l_3^{(2)^{'}}$ in coordinates $w^{(2)}$ is given by $\{
w_2^{(2)} = w_3^{(2)} = 0 \}$. Recall, that $l_0^{(1)}$ is given by $\{
z_2^{(1)} = z_3^{(1)} = 0 \}$. So from (1.1.5) we see that $\phi $ sends the
center of blown up $l_0^{(1)}$ into another one $l_3^{(2)^{'}}$, and thus the
lifting $\hat \phi $ is holomorphic.

Now take small $\varepsilon > 0 $. By $X_0^{\pm {\varepsilon }}$ denote the
ball in $\cc^3$  of radii ${1 \over 2}\pm \varepsilon $. By $X_j^{\pm
\varepsilon}, j = 1,2,3 $ denote manifolds obtained by our blowings up from
$X_0^{\pm \varepsilon }$ as $X_j$ from $X_0$. By $X_j^{(i)^{\pm \varepsilon }}$
 denote the corresponding manifolds obtained from $X_0^{(i)}$, i = 1,2.

Take $ U^{\varepsilon } = \hat X_3^{(2)^{\varepsilon }}\setminus
\hat \phi (\hat X_0^{(1)^{- \varepsilon }})$ - domain in
$\hat X_3^{{(2)}^{\varepsilon }}$. Here $(2)^{\varepsilon }$ and
$(1)^{- \varepsilon }$ one should understand as indices, not as numbers. By
$\pi_0^{(1)} : \hat X_0^{(1)} \longrightarrow X_0^{(1)}$, $\pi _i^{(2)} :
X_{i+1}^{(2)} \longrightarrow X_i^{(2)}, i = 0,1,2, \pi^{(2)} =
\pi_0^{(2)}\circ \pi_1^{(2)}\circ \pi_2^{(2)}$ we denote the corresponding
modifications.

Consider a map $\psi : X_0^{(2)^{\varepsilon }} \longrightarrow
X_0^{(1)^{\varepsilon }}$ which in our coordinates $z^{(i)}$ is given as
identity:

$$z_1^{(1)} = z_1^{(2)}$$
$$z_2^{(1)} = z_2^{(2)}     \eqno(1.1.6)$$
$$z_3^{(1)} = z_3^{(3)}$$

Lift this map onto the blowings up to get

$$
\hat \psi :  \hat X_0^{(2)^{\varepsilon }}\longrightarrow
\hat X_0^{(1)^{\varepsilon}} \eqno(1.1.7)
$$

Now remark that further blowings up of $\hat X_0^{(2)^{\varepsilon }} =
X_1^{(2)^{\varepsilon }}$ do not effect the neighborhoods of the boundaries of
$\hat X_0^{(i)}$. We shall denote this neighborhoods as $V^{(i)^{\varepsilon }}
$, having in mind that they are copies in $\hat X_0^{(i)^{\varepsilon }}$ of
 $V^{\varepsilon} = \hat X_0^{\varepsilon } \setminus cl(\hat X_0^{-
\varepsilon })$. So biholomorphism $\hat \psi $, defined in
(1.1.7),can be restricted to a biholomorphism (which we denote by the same letter)
$\hat \psi : V^{(2)^{\varepsilon }} \longrightarrow V^{(1)^{\varepsilon }}$
. Now we can define a biholomorphism $g = \hat \phi \circ \hat \psi $ between
$V^{(2)^{\varepsilon }}$ and $\hat \phi (V^{(1)^{\varepsilon }})$ -
neighborhoods of two connected components of the boundary of $U^{(2)}$ in
$U^{(2)^{\varepsilon }}$. Here we defined $U^{(i)}$ and $U^{(i)^{\varepsilon }
}$ to be a copies in $X_3^{(i)^{\varepsilon }}$, of the open sets $ U =
\hat X_3^{(2)} \setminus \hat \phi (\hat X_0^{(1)})$ and $U^{\varepsilon } =
\hat X_3^{(2)^{\varepsilon }}\setminus \hat \phi
(\hat X_0^{(1)^{- \varepsilon }})$.

Now we can glue the neighborhoods $V^{(2)^{\varepsilon }}$ and $\hat \phi
(V^{(1)^{\varepsilon }})$ by biholomorphism $g$ to get from $U^{(2)^{
\varepsilon }}$ a compact complex threefold $X$.

Note that after gluing $E_0$ and $E_3$ constitute one divisor. We shall
denote it both as $E_0$ or $E_3$.
\bigskip
\noindent
{\sl 1.2. Construction of the universal cover and nonextendable map.}
\smallskip
\rm
We must prove two facts about $X$. First: that there exists a meromorphic
mapping $F : B^3 \setminus \{ 0\} \longrightarrow X$ which does not extend
to the origin. Second: any meromorphic map $f : H^2(r) \longrightarrow X$
extends to a meromorphic map $\hat f : \Delta ^2 \longrightarrow X$. Here
$H^2(r) = \{ (u_1,u_2) \in \Delta ^2 : \vert u_1\vert < r$  or $1-r < \vert
u_2\vert < 1 \}$ is a Hartogs figure of dimension two. To prove this we shall
need the universal covering $\widetilde{X}$ of $X$. In this paragraph we give
the construction of $\widetilde{X}$ together with nonextendable mapping
$F : B^3 \setminus \{ 0\}  \longrightarrow X$.

We shall make this construction in 5 steps.

In the previous paragraph we had constructed a biholomorphic mapping
$g : V^{\varepsilon } \longrightarrow \hat \phi (V^{\varepsilon })$, where
both $V^{\varepsilon }$ and $\hat \phi (V^{\varepsilon })$ are open subsets of
$\hat X_3^{\varepsilon }$, and actually are the neighborhoods of the connected
components of the boundary of $U = \hat X_3 \setminus \hat \phi (\hat X_0)$.
To simplify a bit our notations we  denote the subset $\hat \phi
(V^{\varepsilon })$ of $\hat X_3^{\varepsilon }$ simply as $V_0^{\varepsilon}$.
\smallskip
\noindent
Step 1. Construction of manifolds $X^{(i)}$ and $U_i^{\varepsilon }, i\ge 0.$

Put $X^{(0)} = \hat X_3^{\varepsilon }$. Domain $U^{\varepsilon }$ in
$X^{(0)}$ denote as $U_0^{\varepsilon }$. Attach $X^{(0)}$ to
$U^{\varepsilon }$ by $g : V^{\varepsilon } \longrightarrow
V_0^{\varepsilon }\subset U^{\varepsilon }$. Manifold obtained, denote as
$X^{(1)}$.

Manifold $X^{(1)}$ has exactly one end, namely $V^{\varepsilon }$. Domain in
$X^{(1)}$ which is $U^{\varepsilon }\sqcup U_0^{\varepsilon }/g$ we denote by
$U_1^{\varepsilon }$.

We can repeat this construction with $X^{(1)}$ instead of $X^{(0)}$. Namely,
take $X^{(1)}$ and attach it to $U^{\varepsilon }$ by $g : V^{\varepsilon }
\longrightarrow V_0^{\varepsilon } \subset U^{\varepsilon }$ to get $X^{(2)}$.
 Domain in $X^{(2)}$ which is $U_1^{\varepsilon }\sqcup U^{\varepsilon }/g$
we denote by $U_2^{\varepsilon }$.

Suppose that the manifold $X^{(i)}$ is constructed and it has one end
$V^{\varepsilon }$. Attach $X^{(i)}$ to $U^{\varepsilon }$ by $g : V^
{\varepsilon }\longrightarrow U^{\varepsilon }$ to get $X^{(i+1)}$, which
again has one end $V^{\varepsilon }$. Domain in $X^{(i+1)}$ which is
$U_i^{\varepsilon }\sqcup U^{\varepsilon }/g$ denote by $U_{i+1}^{\varepsilon}$
.
\smallskip
\noindent
Step 2. Construction of manifolds $X^{(-i)}$ and $U_{-i}^{\varepsilon },i\ge 0
$.

As $X^{(0)}$ we take again $\hat X_3^{\varepsilon }$ with $U_0^{\varepsilon }
= U^{\varepsilon }$.

To obtain $X^{(-1)}$ attach $\hat X_3^{\varepsilon }$ to $U_0^{\varepsilon }$
by $g : V^{\varepsilon }\longrightarrow V_0^{\varepsilon }\subset
U_0^{\varepsilon}$. Domain in $X^{(-1)}$ which is $U_0^{\varepsilon }\sqcup
U^{\varepsilon }/g$ denote by $U_{-1}^{\varepsilon }$.

Suppose that the manifold $X^{(-i)}$ and domain $U_{-i}$ in $X^{(-i)}$ are
constructed. To get $X^{(-i-1)}$ attach $\hat X_3^{\varepsilon }$ to $U_{-i}^
{\varepsilon }$ by $g : V^{\varepsilon }\longrightarrow V_0^{\varepsilon }
\subset U_{-i}^{\varepsilon }$ and put $U_{-i-1}^{\varepsilon } = U_{-i}^
{\varepsilon } \sqcup U^{\varepsilon }/g$.
\smallskip
\noindent
Step 3. Construction of the universal cover.

Manifolds $U_i^{\varepsilon }$ and $U_{-i}^{\varepsilon }$ have the common
part $U_0^{\varepsilon }$. So we can glue them by the identification in
$U_i^{\varepsilon }\sqcup U_{-i}^{\varepsilon }$ of two copies of $U_0
^{\varepsilon }$. The manifold obtained we denote by $U_{-i,i}^{\varepsilon }
$. In the same way we can consider $U_i^{\varepsilon }\sqcup U_{-j}^
{\varepsilon }/g = U_{i,-j}^{\varepsilon }$.

We have now an increasing sequence of complex manifolds  $U_0^{\varepsilon }
\subset \subset U_{1,-1}^{\varepsilon } \subset  \subset ...  ...
 \subset \subset U_{i,-i}^{\varepsilon }   \subset  \subset \cdot \cdot
\cdot $ . The union $\bigcup_{i\ge 0}U_{i,-i}^{\varepsilon } = U_{\infty }$
is a complex manifold, and we are going to
prove that it is the universal cover of $X$. We first remark, that obviously
$U_{\infty }$ is simply connected.

\smallskip
\noindent
Step 4. Construction of group of covering transformations.

We start with the biholomorphic mapping $g : V^{\varepsilon }\longrightarrow
V_0^{\varepsilon }$ which sends one end of $U_0^{\varepsilon }$ to another.
Recall that $U_0^{\varepsilon }/g \cong X$ is our compact threefold.
\smallskip
\noindent
{\bf Lemma 1.2.1} $g$ extends to biholomorphic automorphism of $U_{\infty }$.
Moreover we have \hfill\break $X \cong  \- U_{\infty }/<g^n>_{n\in Z}$.

\noindent
Proof. The heuristic argument for the extendability of $g$ is that it is given
"as identity" in globally defined coordinates. But we shall use here the
Hartogs-type extension argument. After this the second statement of the lemma
becomes clear.

For any pair of integers $(i,j)$ such that $i > j$ we consider a subdomain
$U_{i,j}^{\varepsilon }$ in $U_{\infty }$ defined in a following way.

If $i\ge 0\ge j$ then $U_{i,j}^{\varepsilon }$ is already defined.

If $0\ge i > j$ then $U_{i,j}^{\varepsilon } = U_j^{\varepsilon }\setminus
U_i^{{\varepsilon }-}$.

If $i > j \ge 0$ then $U_{i,j}^{\varepsilon } = U_i^{\varepsilon }\setminus
U_j^{-{\varepsilon }}$.

Here domains $U_i^{{\varepsilon }-}$ are constructed from $U^{{\varepsilon }-}
 = \hat X_3^{(2)^{\varepsilon }} \setminus \hat \phi (\hat X_0^{(1)^{
\varepsilon }})$ in a similar way as $U_i^{\varepsilon }$ from $U^
{\varepsilon }$, and $U_j^{-{\varepsilon }}$ form $U^{-{\varepsilon }} =
\hat X_3^{(2)^{-{\varepsilon }}} \setminus \hat \phi (\hat X_0^{(1)^
{-{\varepsilon }}})$.

By attaching $X_3^{\varepsilon }$ to $U_{i,j}^{\varepsilon }$ by $g$ we get
manifolds $X^{(i,j)}$, $\infty \ge i > j\ge -\infty $.

Our map $g : V^{\varepsilon }\longrightarrow V_0^{\varepsilon }$ we consider
as a map from the open subset of $U_{0,-\infty }^{\varepsilon }$ into itself.
Note that the image of $g$ is contained in $U_{-1,-\infty}^{\varepsilon }$. We
are going to prove that $g$ extends to a biholomorphic map from $U_{0,-\infty }
^{\varepsilon }$ onto $U_{-1,-\infty }^{\varepsilon }$.

For this purpose take any $j < 0$. Remark that $g$ maps the neighborhood $V^{
\varepsilon }$ of the boundary of $X^{(0,j)}$ biholomorphically onto the
neighborhood $V_0^{\varepsilon }$ of the boundary of $X^{(-1,j-1)}$. Note that
$X^{(0,j)}$ and $X^{(-1,j-1)}$ could be naturally blown down to get an $X_0$
 - the usual ball in $\cc^3$. So $g$ became a biholomorphism of the neighborhood
 of $S^5$ onto itself. So it extends by the Hartogs extension theorem.

In the same way $g^{-1}$ extends to a biholomorphic map of $U_{\infty,-1}^
{\varepsilon }$ onto $U_{\infty,0}^{\varepsilon }$.

\smallskip
\hfill{q.e.d.}

By $p : U_{\infty }\longrightarrow X$ we denote the covering map.

\smallskip
\noindent
Step 5. Construction of nonextendable map.

Our construction gives a natural bimeromorphic map $\widetilde{f}$ from
$X_0^{\varepsilon }\setminus \{ 0\}\cong \bb^3_*$ onto $U_{0,-\infty }^
{\varepsilon }$.
Then taking the composition $p\circ \widetilde{f} = F$ we get our nonextendable
 map.

Another way to obtain $\widetilde{f} $ is to consider the composition $i\circ
\pi ^{-1}$ of inclusion $i : V^{\varepsilon }\longrightarrow U_{0,\infty }^
{\varepsilon }$ with the restriction of blowing up $\pi ^{-1}$ onto $X_0^
{\varepsilon }\setminus cl(X_0^{-\varepsilon })$. Now, similarly to the step 4,
one can easily observe that $i\circ \pi ^{-1}$ extends to a bimeromorphic
map $\widetilde{f} : X_0^{\varepsilon }\setminus \{ 0\} \longrightarrow U_
{0,-\infty }^{
\varepsilon }$.

\bigskip
\noindent
{\sl 1.3. Two - dimensional extension property of $X$.}
\smallskip

In this paragraph we shall prove that every meromorphic map $f : H^2(r)
\longrightarrow X$ extends to a meromorphic map $\hat f : \Delta ^2
\longrightarrow X$.

Take $0 < \delta < r$ and denote by $H_{\delta }(r) = \{ (u_1,u_2)\in
\Delta ^2 :\vert u_1\vert < r - \delta , \vert u_2\vert < 1 - \delta $ or $
\vert  u_1\vert < 1 - \delta , 1 - r + \delta < \vert u_2\vert < 1 -
\delta \}$ - shrunked Hartogs domain. It is sufficient to prove the
extendability of $f$ onto $\Delta _{1-\delta }^2 = \{ (u_1,u_2) \in C^2 :
\vert u_j\vert < 1 - \delta ,j = 1,2\}$ for arbitrary $\delta > 0$.

Using the fact that $H^2(r)$ is simply connected we can lift $f$ to
a meromorphic map $\widetilde{f} : H^2(r)\longrightarrow U_{\infty }$.
$H_{\delta }(r)\subset \subset H^2(r) $ so $\widetilde{f} (H_{\delta }
(r)) \subset \subset U_{\infty }$, the latter means that there exists $i <
\infty $ such that $\widetilde{f} (H_{\delta }\}(r)) \subset U_{i,-i}^
{\varepsilon }$. Applying $g^i$ to $\widetilde{f} $ we get a map $f_1 =
g^i\circ \widetilde{f} : H_{\delta }(r) \longrightarrow U_{0,-2i}^
{\varepsilon }$. If we extend $f_1$ as a meromorphic map $\hat f_1$ from
$\Delta _{1-\delta }^2$ to $U_{0,-\infty }^{\varepsilon }$ then $p\circ
g^{-i}\circ \hat f_1 $ will give us a desired extension of $f$.

Take a composition $h = p \circ f_1 : H_{\delta}(r) \longrightarrow X_0 \cong
 B_{1/2}^3 \subset \cc^3$. This is a holomorphic mapping of $H_{\delta }(r)$
into  $X_0 \subset \cc^3$. So it extends to a holomorphic map $\hat h : \Delta
_{1-\delta }^2 \longrightarrow X_0$. By $A$ denote the preimage of zero under
$\hat h$. There are three possibilities.

\noindent
{\bf Case 1.} $A = \Delta _{1-\delta }^2$ . That means that $f(H_{\delta }(r))
\subset E$, where $E$ is an exceptional divisor of $X$. But all three
components of $E$ are bimeromorphic to $\cc\pp^2$. So in that case $f$ obviously
extends onto $\Delta _{1-\delta }^2$.

\noindent
{\bf Case 2.} $A$ has components of positive dimension.

Denote them by $A_1,...,A_N$. Let also $a_1,...,a_M$ are the isolated preimages
of zero under $\hat h$.

Because $p : U_{0,-\infty }^{\varepsilon }\setminus \biggl[\bigcup_{i=1}^
{-\infty }
\bigcup_{j=0}^3 E_j^i\cup \bigcup_{j=1}^3E^0_j\biggr] \longrightarrow
X_0\setminus \{ 0\}$ is bimeromorphic
we have an extension of $f_1$ onto $\Delta _{1-\delta }^2\setminus \biggl[
\bigcup_{j=1}^NA_j^{-}\cup \bigcup_{k=1}^M\{ a_k\} \biggr]$. Here we denote
by $A_j^{+} = A_j \cap H_{\delta }(r)$, by $A_j^{-} = A_j\setminus A_j^{+}$
,and  $E_j^{i}$ denote the copies of $E_j$ in $U_{i,i-1}^{\varepsilon }, j =
0,1,2,3$. Note that $E_3^i$ is glued to $E_0^{i-1}$ and constitute a single
divisor, which we shall denote either by $E_3^i$ or $E_0^{i-1}$.

Consider now $f_1$ as a meromorphic map from $H_{\delta }(r)$ to $U_{0,-2i}^
{\varepsilon } \subset X^{(-2i-L)}$. Because $X^{(-2i-L)}$ is nothing but a
blown up ball we can extend $f_1$ as a map $\hat f$ from $\Delta _{1-\delta }^
2$ to $X^{(-2i-L)}$. Now $\hat f_1(\bigcup_{j=1}^NA_j) \subset \bigcup_{k=0}^
{-2i}\bigcup_{j=0}^3E_j^k$, because $f_1(\bigcup_{j=1}^NA_j^{+}) \subset
\bigcup_{k=0}^{-2i}\bigcup_{j=0}^3E_j^k$. That means that $\hat f_1(\bigcup_
{j=1}^NA_j)$ do not intersect $E_2^{-2i-L}$ for $L$ big enough. This follows
from Lemma 1.3.1 stated below and the fact that nontrivial holomorphic
mapping from 1-dimensional space is always locally proper. So its image is
locally contained in a curve by the Remmert proper mapping theorem. Before
applying Lemma 1.3.1 one only should embed this curve into a (nonsmooth in
general) hypersurface. So if we consider
$f_1$ as a map into $U_{0,-\infty }$, then $f_1(\bigcup_{j=1}^NA_j) \subset
U_{0,-2i-2}^{\varepsilon }$. So there is a neighborhood $W$ of $\bigcup_{j=1}^
NA_j$ such that $f_1(W) \subset U_{0,-2i-2}^{\varepsilon }$. The latter means
that $p\circ g^{-i}\circ f_1 $ extends onto $W$. So our map $f$ is extended
onto $\Delta _{1-\delta }^2 \setminus \bigcup_{k=1}^M\{a_k \}$.

So the last case left.

\noindent
{\bf Case 3.} $A$ is a finite set of points $\{a_1,...,a_M\}$.

Take a neighborhood $B$ of $a_1$. Because $a_1$ is isolated preimage of zero
under $\hat h$,there is a neighborhood of $a_1$,say $B$ itself, such that $\hat
 h\mid _{\bar B}$ is proper. So $\hat h(\bar B)$ is contained in a germ of
hypersurface $P$ in $X_0$, passing through origin. Note that this fact was not
obvious from the beginning because of the well known Osgood example, see [4],
p.155. If we shall prove that the lifting $p ^{-1}(P)$ of $P$ onto
$U_{0,-\infty }$
 is relatively compact in $U_{0,-\infty }^{\varepsilon }$ then we are done,
because then $f_1$ maps $B\setminus \{a_1\}$ into some $U_{0,-L}^{\varepsilon }
 \subset X_{1,-L}$.

To prove that the lifting of $P$ by $p$ is relatively compact in $U_{0,-
\infty }^{\varepsilon }$, the next lemma is sufficient. Denote by $\Phi $ the
transformation given by (1.1.4).

\noindent
{\bf Lemma 1.3.1} \it Let $P(z)$ be a convergent series in variables $z_1,
z_2,z_3$.
 Then ,there exist an $N$ such that $P(\Phi ^N(z))$ has the form $z_1^{N_1}
\cdot z_2^{N_2}\cdot z_3^{N_3}\cdot \biggl[ a + F(z)\biggr]$, where $a$ is some nonzero
constant, $F$ is convergent series and $F(0) = 0$.

\smallskip\noindent\sl
Proof. \rm Consider a Newton polyhedron $N_P$ of $P(z)$ in $\rr^3$. Consider
also
the next unimodular transformation
$$A =
\left(\matrix{
1 & 1 & 0 \cr
1 & 2 & 1 \cr
1 & 1 & 1 \cr
}
\right)
$$

\noindent
of $\rr^3$. One can easily check the following properties of $A$ :

1. The real eigenvalue $\mu $ of $A$ satisfies $\mu > 3$.

2. Consequently the complex one $\lambda $ satisfies $\vert \lambda \vert ^2 =
1/{\mu } < 1/3$.

3. The eigenvector $v$, which corresponds to $\mu $, has positive coordinates.

4. The $A$ - invariant plane $\Pi $, which corresponds to $\lambda $ intersects
the positive octant by zero.

5. $\Pi $ contains no vectors with integer coordinates.

The latter is true because $A$ acts on $\Pi $ as rotation together with
contraction by $\lambda $.

From 5 we see that there is a constant $\alpha > 0$, such that every edge of
$N_P$ has an angle at least $\alpha $with $\Pi $. Now, because $A$ acts as
contraction
on $\Pi $ and as dilatation on [$Rv$],all edges of $N_{P^N}$ have angles with
$\Pi $ which are uniformly tending to the angle between $v$ and $\Pi $. Here
by $P^N$ we denoted $P(\underbrace{\Phi \circ ...\circ \Phi }_{N}(z))$.

This proves the lemma.

\hfill{q.e.d.}
\smallskip
\noindent
{\bf Remarks} 1. The example, just constructed, has some common futures with
an example of Hirschowitz, see [Hs], of weakly but not strongly meromorphic
mapping.
 In fact the map $F$, obtained in $1.2$ is weakly meromorphic in the sense,
that for any complex curve $C\ni 0$ the restriction $F\mid _{C\setminus 0}$
extends holomorphically onto $C$. But clearly the two dimensional extension
property of our $X$ is much stronger condition than just a weak
meromorphicity of some given map, namely of $F$, which is not strongly
meromorphic.

\smallskip\noindent 2.
Note also, that $X$ contains some sort of "spherical shell". Namely a \it
bimeromorphic \rm image $V^{\varepsilon }$ of the neighborhood of $\ss^5$ from
$\cc^3$ which doesn't bound a domain in $X$. We want to point out that this
is different from
the situation examined by M. Kato in [Ka-1], who studied so called "global
spherical shells" i.e. \it biholomorphic \rm images of some neighborhoods of
$S^{2n-1}$ from $\cc^n$, which doesn't bound. Our $X$ doesn't contain
a GSS in the sense of Kato.

\bigskip\noindent\bf
2. Meromorphic polydisks.
\smallskip\noindent\sl
2.1. Generalities on pluripotential theory. \rm

We start with some well known facts from pluripotential theory . Let $D$ be
an open subset of $\cc^n$ and $S$ subset of $D$.
Consider the next class of functions
\smallskip
$${\cal U}(S,D) = \{ u\in {\cal P}_{-}(D) : u|_S\ge 1\}\eqno
(2.1.1)
$$
\smallskip
\noindent
where by ${\cal P}_{-}(D)$ we denote the class of nonnegative
plurisuperharmonic functions in $D$.
\smallskip
\noindent
\bf Definition 2.1.1. \rm The  lower regularization $w_*$ of the function
\smallskip
$$
w(\zeta ,S,D) = \inf \{u(\zeta ) :
u\in {\cal U}(S,D)\}\eqno(2.1.2),
$$
is called a ${\cal P}$- measure of $S$ in $D$, i.e.
$$
w_*(z,S,D)=\lim_{\zeta \rightarrow z}\inf w(\zeta ,S,D) \eqno(2.1.3)
$$
\smallskip\noindent
Note that $w_{\ast }$ is plurisuperharmonic in $D$ .
\smallskip
\noindent
\bf Definition 2.1.2. \rm A point $s_0\in S$ is called a locally regular point
of $S$ if $w_*(s_0,S\cap \Delta (s_0,\varepsilon ),\Delta (s_0,\varepsilon ))
= 1$ for all $\varepsilon >0$.

We shall also say that the set $S$ is locally regular at $s_0$.

The next statement is known as the two constants theorem; see for ex. [Sa]:

\it let $v$ be plurisubharmonic in $D$ such that $v|_S\le M_0$ and
$v|_D\le M_1$.
Then for $z\in D$ one has
\smallskip
$$
v\le M_0\cdot w_*(z) + M_1\cdot [1-w_*(z)]\eqno(2.1.4).
$$
\smallskip
\rm The following lemma can be proved by Taylor expansion and using (2.1.4).
\smallskip
\noindent
\bf Lemma 2.1.1. \it Suppose that the function $f$ is holomorphic and bounded
by modulus with constant $M$ in $\Delta^n \times \Omega $, where $\Omega $ is a
subdomain in $\Delta ^q$. Let for $s\in S$ $ f_s=f(s,\cdot )$ extends
holomorphically onto $\Delta _s^q = \{ s\} \times \Delta ^q$. Suppose that
all this
extensions are also bounded by modulus with $M$. If $s_0\in S$ is a locally
regular point then for every $0<b<1$ there is an $r>0$ such that $f$
holomorphically extends onto $\Delta^n (s_0,r)\times \Delta ^q(b)$ and is
bounded by modulus with $M\cdot {1+b\over 2}$.

\rm
\smallskip
Recall that a compact subset $S\subset D $ is called pluripolar
if there exists a plurisubharmonic function $u$ in $D ,u\not\equiv -\infty $,
such
that $u\mid_S \equiv -\infty $. $S$ is complete pluripolar if $S=\{ z: u(z)=
-\infty \} $. We shall repeatedly use the following statement:

\it if subset $S\subset D $ ($D$ is now pseudoconvex) is not locally regular
at all its points then $S$

is pluripolar,

\noindent
\rm see [B-T],[Sa]. We shall use  the following immediate corollary from the
famous Josefson theorem, see [Kl]:
\smallskip\noindent\bf
\it Let $\Omega $ be a pseudoconvex set in $\cc^n $, and let $S_n$ be a
sequence of subsets of $\Omega $ such that:

1) $S_1$ is closed and pluripolar in $\Omega $;

2) $S_{n+1}\subset \Omega\setminus S_n$ is closed in $\Omega \setminus S_n$
pluripolar;

\noindent
Then $S:=\bigcup_{n=1}^{\infty }S_n$ is pluripolar in $\Omega $.
\smallskip\noindent
\rm

\bigskip
\noindent {\sl 2.2. Passing to higher dimensions: holomorphic case.}
\smallskip\rm
We want to show now that if a complex space $X$ possesses a \it holomorphic
\rm extension property in dimension two than  $X$ possesses this property in
all higher dimensions. Let $n\ge 3$. Put
$$
E^n(r) = (\Delta^{n-2}_r\times \Delta_r\times \Delta )\cup (\Delta^{n-2}_r
\times \Delta \times A_{1-r,1}) = \Delta^{n-2}_r\times H^2(r)
$$
\smallskip\noindent\bf
Lemma 2.2.1. \it If any meromorphic (holomorphic) map $f:E^n(r)\to X$
extends to a meromorphic (holomorphic) map $\hat f:\Delta^{n-2}_r\times
\Delta^2\to X $ then the space $X$ possesses a meromorphic (holomorphic)
extension property in dimension $n$.
\smallskip\noindent\sl
Proof. \rm Let $f:H^n(r)\to X$ be a meromorphic (holomorphic) map. Because
$E^n(r)\subset H^n(r)$ the map $f$ can be extended onto $\Delta^{n-2}_r\times
\Delta^2$. We shall prove by induction on $j=2,\ldots,n$ that $f$ can be extended
onto $\Delta^{n-j}_r\times \Delta^j$. For $j=n$ we shall get the statement of
{\sl Lemma }.

Suppose it is proved for $j$. Let $z_1,\ldots,z_n$ be the coordinates in
$\cc^n$. For a point $z'_0=(z_0^{n-j+1},\ldots,z_0^{n-1})\in \Delta^{j-1}_
{1-r}$
consider a domain
$$
E_{z_0'}^n(r)=[\Delta_r^{n-j-1}\times \Delta_r\times \Delta_r^{j-1}(z_0')
\times
\Delta ]\cup [\Delta_r^{n-j-1}\times \Delta \times \Delta_r^{j-1}(z_0')\times
A_{1-r,1}].
$$
\noindent
Here by $\Delta_r^{j-1}(z_0')$ we denote the $(j-1)$-disk of radii $r$
centered at $z_0'$. Because $E^n_{z_0'}\subset (\Delta_r^{n-j}\times
\Delta^j)\cap H^n(r)$ the map $f$ by the induction hypothesis extends onto
$\Delta_r^{n-j-1}\times \Delta \times \Delta_r^{j-1}(z_0')\times \Delta $.
But
$$
\bigcup_{z_0'\in \Delta^{j-1}}\Delta_r^{n-j-1}\times \Delta \times
 \Delta_{1-r}^{j-1}(z_0')\times \Delta = \Delta_r^{n-j-1}\times \Delta^{j+1}
$$
\noindent and the {\sl Lemma } is proved.
\smallskip
\hfill{q.e.d.}
\smallskip

Now left to prove the following
\smallskip\noindent\bf
Lemma 2.2.2. \it Let $f:E^{n+1}(r)\to X$ be a holomorphic mapping into a
complex space $X$. Suppose that for every $z_0^1\in \Delta_r$ the restriction
$f_{z_0^1}$ of $f$ onto $\{ z^1=z_0^1\} \cap E^{n+1}(r):=E^n_{z_0^1}(r)$
extends
holomorphically onto $\Delta^n_{z_0^1}(r):=\{ z_0^1\} \times \Delta^{n-1}_
r\times \Delta^2 $. Then $f$ holomorphically extends onto $\Delta^{n-1}_r
\times \Delta^2$.
\smallskip\noindent\sl
Proof. \rm Shrinking disks and annuli involved in the definitions of
$E^{n+1}(r)$ and $E^n_{z_0^1}(r)$ we can suppose that $f$ is defined in the
neighborhood of $\overline{E^{n+1}(r)}$ and for every $z_0^1\in \bar \Delta_r$
extends holomorphically to the neighborhood of $\{ z_0^1\} \times \overline
{\Delta_r^{n-2}\times \Delta^2} $ on the hyperplane $\{ z^1=z_0^1\} $, which
doesn't
depend on $z_0^1$. Making homothety by ${1\over r}$ on first $(n-1)$
coordinates we can suppose that $f$ is holomorphic in the neighborhood of the
closure of domain $(\Delta^{n-1}\times \Delta_r\times \Delta )\cup (\Delta^
{n-1}\times \Delta \times A_{1-r,1})$ and for every $z_0^1\in \bar \Delta $
the restriction $f_{z_0^1}$ extends holomorphically to the neighborhood of the
closure of $\{ z_0^1\} \times \Delta^{n-2}\times \Delta^2 $ on the hyperplane
$\{ z^1=z^1_0\} $ not depending on $z_0^1$.

It is enough to prove that for every $z_0^1\in \bar\Delta $ the map $f$
extends
holomorphically onto the $(n+1)$-dimensional neighborhood of $\{ z_0^1\} \times
 \bar\Delta^n:=\bar\Delta^n_{z_0^1}$. Because $f_{z_0^1}$ is holomorphic on
$\bar\Delta^n_{z_0^1}$ the graph $\Gamma_{f_{z_0^1}}$ is an imbedded
$n$-disk in the space $\Delta^{n+1}\times X$. Let $W$ be a Stein neighborhood
of $\Gamma_{f_{z_0^1}}$, see [Si-3]. From the holomorphicity of $f$ on
 $(\Delta^{n-1}\times \Delta_r\times \Delta )\cup (\Delta^
{n-1}\times \Delta \times A_{1-r,1})$ follows that for $z^1$ in some
neighborhood $U\ni z_0^1$, $f_{z^1}((\Delta^{n-1}\times \Delta_r\times \Delta
)\cup (\Delta^{n-1}\times \Delta \times A_{1-r,1}))\subset W$. But $W$ is
Stein, so
$f_{z^1}(\bar\Delta^{n-2}\times \Delta^2)\subset W$ for all $z_0^1\in U $. I.e.
 the map $f$ maps $U\times \bar\Delta^{n-2}\times \bar\Delta^2$ to $W$. So the
 statement of the {\sl Lemma } is reduced to the analogous statement for
mappings into Stein manifolds, and for them our {\sl Lemma} is obvious.
\smallskip
\hfill{q.e.d.}

\smallskip
As we already had mentioned in the {\sl Introduction}, our example shows
that the analogous statement
for {\it meromorphic} mappings is not valid.
Here we only want to point out that the proof of {\sl Lemma 2.2.2} for the
{\it meromorphic} case fails, because if $f_{z_0^1}$ is essentially meromorphic
 (i.e. not holomorphic) then $\Gamma_{f_{z_0^1}}$ doesn't have Stein
neighborhoods.

\smallskip\noindent
\sl  2.3. Sequences of analytic sets of bounded volume.\rm
\medskip\rm
\smallskip\rm
A Hermitian metric form on the complex space $X$ we define in the following
way. Let an open covering $U_{\alpha }$ of $X$ is given together with proper
holomorphic injections $\phi_{\alpha }:U_{\alpha }\to V_{\alpha }$ into a
domains $V_{\alpha }\subset \cc^{n(\alpha )}$. Let $U'_\alpha$ be the
images of $U_{\alpha }$. Let $\{ w_{\alpha }\} $ are positive (1,1)-forms on
$V_{\alpha }$. $\{ w_{\alpha }\} $ defines a Hermitian metric form on $X$ if
 $(\phi_{\alpha }\circ \phi_{\beta }^{-1})^*w_{\beta }=w_{\alpha }$ for all
$\alpha ,\beta $. Note that  $\phi_{\alpha }\circ \phi_{\beta }^{-1}$ is
defined in some neighborhood of $\phi_{\beta }(U_{\alpha }\cap U_{\beta })$
in $\cc^{n(\beta )}$. We say that the metric $w$ is K\"ahler if $dw_{\alpha
 }=0$ on $V_{\alpha }$ for all $\alpha $.
\smallskip
Fix a complex space $X$, equipped with some Hermitian metric $h$.
By $w_h$, or
simply by $w$ denote the $(1,1)$-form canonically associated with $h$.
 Let $\Delta ^q$ be a polydisk in $\cc^q$ with standard
Euclidean metric $e$. The associated form will be denoted by $w_e = dd^c
\Vert z\Vert ^2 = i/2 \sum_{j=1}^qdz_j\wedge d\bar z_j $. By $p_1:\Delta ^q
\times X \longrightarrow \Delta ^q$ and $p_2:\Delta ^q\times X\longrightarrow
X$ we denote the projections onto the first and second factors. On the product
$\Delta \times X$ we consider the metric form $w = p^{\ast }_1w_e + p^{\ast }_2
w_h$.
\smallskip
\noindent
\bf Definition 2.3.1. \rm  By a meromorphic $q$-disk in the complex space $X$
we shall understand a meromorphic mapping $\phi :\Delta ^q \longrightarrow X $
, which is defined in some neighbourhood of the closure $\bar \Delta ^q$.
\smallskip
It will be convenient for us to consider instead of mappings $\phi : \Delta ^q
\longrightarrow X$ their graphs $\Gamma _{\phi }$. By $\hat \phi = (z,\phi
(z))$ we shall denote the mapping into the graph $\Gamma _{\phi } \subset
\Delta ^q\times X$. The volume of the graph $\Gamma _{\phi }$ of the mapping
$\phi $ is given by
\smallskip
$$
{\sl vol}(\Gamma _{\phi }) = \int_{\Gamma _{\phi }}w^q =
\int_{\Delta ^q}(\phi ^{\ast }w_h + dd^c\Vert z\Vert ^2)^q \eqno(2.3.1)
$$
\smallskip
\noindent
Here by $\phi ^{\ast }w_h$ we denote the preimage of $w_h$ under $\phi $, i.e.
$\phi ^{\ast }w_h = (p_1)_{\ast }p^{\ast }_2w_h$.

Recall that the Hausdorff distance between two subsets $A$ and $B$ of the
metric space $(Y,\rho )$ is a number $\rho (A,B) = \inf \{ \varepsilon : A^
{\varepsilon } \supset B, B^{\varepsilon }\supset A\} $. Here by $A^
{\varepsilon }$ we denote the $\varepsilon $-neighborhood of the set $A$, i.e.
$A^{\varepsilon } = \{y\in Y: \rho (y,A) < \varepsilon \}$.

Further, let $\{ \phi _r\} $ be the sequence of meromorphic mappings of the
complex space $D$ into the complex space $X$.
\smallskip
\noindent
\bf Definition 2.3.2. \rm   We shall say that $\{\phi _r\}$ converge on the
compacts
 in $D$ to the meromorphic mapping $\phi :D\longrightarrow X $, if for every
relatively compact open $D_1 \subset \subset D$ the graphs $\Gamma _{\phi _r}
\cap (D_1\times X)$ converge in the Hausdorff metric on $D_1\times X$ to the
graph $\Gamma _{\phi }\cap (D_1\times X)$.

First we shall  prove the following
\smallskip
\noindent
\bf Lemma 2.3.1. \it Let $\{\phi _r\}$ be a sequence of meromorphic $q$-disks
in
complex space $X$. Suppose that there exists a compact $K\subset X$ and a
constant $C<\infty $ such that:

a) $\phi _r(\Delta ^q)\subset K$ for all $r$;

b) ${\sl vol}(\Gamma _{\phi _r})\le C$ for all $r$.
\smallskip
\noindent
Then there exists a subsequence $\{\phi _{r_j}\}$ and a proper analytic set $A
\subset \Delta ^q$ such that:

1) the sequence $\{\Gamma _{\phi _{r_j}}\}$ converges in the Hausdorff metric
to the analytic subset $\Gamma $ of $\Delta ^q\times X$ of pure dimension $q$;

2) $\Gamma = \Gamma _{\phi }\cup \hat \Gamma $,where $\Gamma _{\phi }$ is the
graph of some meromorphic mapping $\phi :\Delta ^q \longrightarrow X$,and
$\hat \Gamma $ is a pure $q$-dimensional analytic subset of $\Delta ^q\times X$
,mapped by the projection $p_1$ onto $A$;

3) $\phi _{r_j}\longrightarrow \phi $ on compacts in $\Delta ^q\setminus A$;

4) one has
$$
\lim_{j \longrightarrow \infty }{\sl vol}(\Gamma _{\phi _{r_j}})\ge
{\sl vol}(\Gamma _{\phi }) + {\sl vol}(\hat \Gamma ).\eqno(2.3.2)
$$

5) For every $1\le p\le dimX - 1$ there exists a positive constant $\nu _p
= \nu _p(K,h)$ such that the volume of every pure $p$-dimensional compact
analytic subset of $X$ which is contained in $K$ is not less then $\nu _p$.

6) Put $\hat \Gamma = \bigcup_{p=0}^{q-1}\Gamma _p $, where $\Gamma _p$ is a
union of all irreducible components of $\hat \Gamma $ such that ${\sl dim}
[p_1(\Gamma _p)] = p$. Then
\smallskip
$$
{\sl vol}_{2q}(\hat \Gamma ) \ge \sum_{p=0}^{q-1}{\sl vol}_{2p}(A_p)\cdot
\nu _{q-p}\eqno(2.3.3)
$$
\smallskip
\noindent
where $A_p = p_1(\Gamma _p)$.
\smallskip
\noindent
\bf Proof. \rm 1) The first statement of this lemma is exactly the theorem of
Bishop about sequences of analytic sets of bounded volume, see [St]. Because
$\Gamma _{\phi _{r_j}}\subset \Delta ^q\times K$, so also a limit $\Gamma
\subset \Delta ^q\times K$. Consequently $p_1\mid _{\Gamma }:\Gamma
\longrightarrow \Delta ^q$ is proper. Further, because $p_1(\Gamma _{\phi
_{r_j}}) = \Delta ^q$ for all so $p_1(\Gamma ) = \Delta ^q$. That is our map
$p_1\mid_{\Gamma }:\Gamma \longrightarrow \Delta ^q$ is surjective.

Let $\Gamma \subset \Delta ^q\times X$ be an analytic subset and let $p_1
\mid_{\Gamma }\Gamma : \longrightarrow \Delta ^q$ be the restriction of the
natural
 projection onto $\Gamma $. We shall say that $p_1\mid_{\Gamma }$ is
one to one over a generic point, if there exists an open $D\subset \Delta ^q$
such that $p_1\mid_{\Gamma \cap p_1^{-1}(D)} : \Gamma\cap
p_1^{-1}(D)\longrightarrow D$ is one
to one. In this case, using the Remmert proper mapping theorem one can show
that $D$ can be taken to be $\Delta ^q$ minus proper analytic subset.

2) Define $\Gamma _1 = \{ (z,x)\in \Gamma :{\sl dim}_{(z,x)}(p_1\mid _{\Gamma }
)^
{-1}(z)\ge 1\}$. $\Gamma _1$ is an analytic subset of $\Gamma $. By the
properness of $p_1\mid _{\Gamma }$ $A_1 = p_1(\hat \Gamma _1)$ is analytic
subset of $\Delta ^q$. By $\hat \Gamma $ denote the union of irreducible
components of $\Gamma _1$ of dimension $q$ ant put $A = p_1(\Gamma )$. Then
$\Gamma = \Gamma _{\phi }\cup \hat \Gamma $, where $\Gamma _{\phi }$ is the
union of those irreducible components of $\Gamma $ which are not in $\hat
\Gamma $. Remark that because $dim A \le q-1$ the restriction $p_1\mid _
{\Gamma _{\phi }}$ is surjective and one to one over the general point. To
prove that $\Gamma _{\phi }$ is a graph of some meromorphic mapping it is
sufficient to show that  $\Gamma _{\phi }$ is irreducible. Suppose that
$\Gamma _2$ is some nontrivial irreducible component of $\Gamma _{\phi }$.
Because $p_1\mid _{\Gamma _{\phi }}$ is one to one over a
general point then $p_1(\Gamma _2) = A_2$ is a proper analytic subset of
$\Delta ^q$. But then $\Gamma _2\subset \hat \Gamma $ by the definition of
$\hat \Gamma $. So $\Gamma _{\phi }$ is irreducible and thus is a graph of a
meromorphic mapping which we denote by $\phi $. Let $\Gamma '$ be the union of
all irreducible components of $\Gamma _1$ of dimension less then $q$. Then
$p_1(\Gamma ') = F$ is exactly the set of points of indeterminacy of the map
$\phi $.

3) $\Gamma_{\phi_{r_j}}\to \Gamma =\Gamma_{\phi }\cup \hat\Gamma $ on
$\Delta^q\times X$. So, because $\hat\Gamma =p_1^{-1}(A)$ we have that
$\Gamma_{\phi_{r_j}}\to \Gamma_{\phi }$ on compacts in $(\Delta^q\setminus
A)\times X$. By definition this means that $\phi_{r_j}\to \phi $ on compacts
in $\Delta^q\setminus A$.

4) This statement follows from the generalisation of the theorem of Bishop
made by Harvey-Shiffman, see [H-S]. The theorem of Harvey-Shiffman states
that the sequence of analytic sets $\{ \Gamma _{\phi _{r_j}}
\}$ , which in our case converges to the analytic set $\Gamma $ in the
Hausdorff metric, converges to $\Gamma $ also in locally flat topology, see
Theorem 3.9 from [H-S]. The latter means, in particular, that for any $(k,k)$-
form $\Omega $ on $\Delta ^q\times X$ with compact support $\int_{\Gamma
_{\phi _{r_j}}}\Omega \longrightarrow \int_{\Gamma }\Omega $ when $j
\longrightarrow \infty $.

If one takes as $\Omega $ the $q^{\rm th}$
 exterior power of $w$ times nonnegative
function $\psi $ with compact support, i.e. $\Omega = \psi \cdot w^q$, and
takes into account that $vol(\Gamma ) = \int_{\Gamma }w^q$, then one gets that
for every compact $K\subset \Delta ^q\times X$:
\smallskip
$$
vol(\Gamma_{\phi _{r_j}}\cap K) = \int_{\Gamma _{\phi _{r_j}}}w^q
\longrightarrow \int_{\Gamma \cap K}w^q =
$$
\smallskip
$$
 = {\sl vol}(\Gamma _{\phi }\cap K) + {\sl vol}(\hat \Gamma \cap K)
\eqno(2.3.4)
$$

Now take into account that $K$ is arbitrary and that the limit in (2.3.4) one
should understand with multiplicities. Withought multiplicities one gets
from (2.3.4) the stated unequality (2.3.2).

5) Proof of this item we shall derive from contradiction. Suppose that there
exists a sequence of pure $p$-dimensional compact analytic subsets $\{ C_j\}
$ in $X$ such that $vol(C_j) \longrightarrow 0$ while $j\longrightarrow
\infty $. Going to a subsequence one obtaines that $C_{j_k} \longrightarrow
C$ in the Haussdorff metrik to the compact $p$-dimensional analytic set $C$.
From (2.3.2) we see that ${\sl vol}(C) \le \lim_{k\rightarrow \infty } \inf
{\sl vol}(C_{j_k}) = 0$
which is
impossible.

6) The relation (2.3.3) reflects the fact that our metric $w$ on the product
is
a sum of the metrics on the factors. Recall that $A=p_1(\hat \Gamma )$, where
$\hat \Gamma $ is a union of the irreducible components of the limit $\Gamma $
other then $\Gamma _{\phi }$. As above denote by $A_p = (p_1|_{\Gamma })(\Gamma
_p)$, where $\Gamma _p$ as above is a union of irreducible components of
$\hat \Gamma $ such that ${\sl dim}A_p = p$. For the proof of (2.3.3) it is
sufficient
to show that
\smallskip
$$
{\sl vol}(\Gamma _p) \ge {\sl vol}(A_p)\cdot \nu _{q-p}\eqno(2.3.5)
$$
\smallskip
Let us underline that $\Gamma _p$ is $q$-dimensional as a union of some $q$-
dimensional components of $\hat \Gamma $, and that $A_p$ is $p$-dimensional.
We have
\smallskip
$$
{\sl vol}(\Gamma _p) = \int_{\Gamma _p}w^q = \int_{\Gamma _p}(p_1^*w_e +
p_2^*w_h)
^q\ge
$$
\smallskip
$$
\ge \int_{\Gamma _p}(p_1^*w_e)^p\wedge (p_2^*w_h)^{q-p} = \int_{A_p}(\int_
{(p_1\mid_{\hat \Gamma })^{-1}(z)}w_h^{q-p})w_e^p \ge \int_{A_p}\nu _{q-p}
\cdot
w_e^p =
$$
\smallskip
$$
= \nu _{q-p}\cdot {\sl vol}(A_p).
$$
\smallskip
In the first unequality we used the fact that all terms of decomposition
$$
(p_1^*w_e + p_2^*w_h)^q = \sum_{j=0}^kC_j^q(p_1^*w_e)^j\wedge (p_2^*w_h)^
{q-j}$$
are positive. In the second one  we used
 that $\nu _{q-p}$ is a minima of volumes of
$(q-p)$-dimensional compact analytic subsets of $X$ contained in $K$.
\hfill{q.e.d.}

\bigskip
\noindent\sl
2.4. Meromorphic families of analytic subsets.
\medskip \rm
Let  $S$ be a set, and $W\subset\subset \cc^q$ an open subset. $W$ is equipped
with the usual Euklidean metric from $\cc^q$.
\smallskip
\noindent
\bf Definition 2.4.1. \it By a family of $q$-dimensional analytic subsets in
complex
space $X$ we shall understand an subset ${\cal F}\subset S\times W\times X $
such that, for every $s\in S$ the set  $ {\cal F}_s = {\cal F}\cap
\{s\}\times
W\times X$ is a graph of a meromorphic mapping of $W$
into $X$.

\smallskip\rm
Suppose further that the set $S$ is equipped with topology and let our space
$X$ be equipped with some Hermitian metric $h$.
\smallskip
\noindent
\bf Definition 2.4.2. \it We shall say that the family ${\cal F}$ is
continuous
at  point $s_0\in S$ if ${\cal H}-\lim_{s\rightarrow s_0}
 {\cal F}_s = {\cal F}_{s_0}$.

\smallskip\rm
Here by ${\cal H}-\lim_{s\rightarrow s_0} {\cal F}_s $
we denote the limit of closed subsets of ${\cal F}_s$ in the Hausdorff
 metric on $W\times X$. ${\cal F}$ is continuous if it is
continuous
at each point of $S$. If $W_0 $ is open in $W$ then the
restriction ${\cal F}_{W_0}$ is naturally defined as ${\cal F}\cap (S
\times W_0 \times X)$.

When $S$ is a complex space itself, we give the following
\smallskip
\noindent
\bf Definition 2.4.3. \it Call the family ${\cal F}$ meromorphic if the
closure $\hat {\cal F}$ of the set
${\cal F}$ is an analytic subset of $S\times W \times X$.

\smallskip\rm
We shall be interested with an interaction of notions of continuity and
meromorphicity of families of meromorphic polydisks.

Let us prove  our main statement about meromorphic families. Consider a
meromorphic mapping $f: V \times W_0 \longrightarrow X$ into a
complex space $X$, where $V$ is a domain in $\cc^p$. Let $S$ be some closed
subset of $V $ and $s_0\in S$
some accumulation point of $S$. Suppose that for each $s\in S$ the
restriction $f_s = f|_{\{ s\} \times W_0} $ meromorphically extends
onto $W\supset\supset W_0$. We suppose additionally that there is a
compact $K\subset\subset X$ such that for all $s\in S$ $f_s(W)\subset K$.

Let, as in Lemma 2.3.1  $\nu _j$ denotes the
minima of
volumes of $j$-dimensional compact analytic subsets contained in our compact
$K\subset X$. Fix some $W_0\subset\subset W_1\subset\subset W$ and put
\smallskip
$$\nu = \min \{ {\sl vol}(A_{q-j})\cdot \nu _j : j = 1,\ldots,q\},\eqno(2.4.1)$$
\smallskip
\noindent
where $A_{q-j}$ are running over all $(q-j)$-dimensional analytic subsets of
$W$, intersecting $\bar W_1$.  Clearly
$\nu > 0$. In the
following {\sl Lemma } the volumes of graphs over $W$
are taken. More precisely, having an Euklidean metric form $w_e=dd^c\Vert
z\Vert^2$ on $W\subset \cc^q$ and Hermitian metric form $w_h$ on $X$, we
consider $\Gamma_{f_s}$ for $s\in S$ as an analytic subsets of $W\times X$
and their volumes are
$$
vol(\Gamma_{f_s}) = \int_{\Gamma_{f_s}}(p_1^*w_e+p_2^*w_h)^q =
\int_W(w_e+(p_1)_*p_2^*w_h)^q,
$$
\noindent
where $p_1:W\times X\to W$ and $p_2:W\times X\to X$ are natural projections.
\smallskip
\noindent
\bf Lemma 2.4.1. \it Suppose  that there exists
a neighbourhood $U\ni s_0$ in $V$ such that, for all $s_1,s_2\in S\cap U$
\smallskip
$$
| {\sl vol}(\Gamma _{f_{s_1}}) - {\sl vol}(\Gamma _{f_{s_2}})
| < \nu /2. \eqno(2.4.2)
$$
\smallskip
\noindent
If  $s_0$ is a locally regular point of $S$ then  there exists
a neighbourhood $V_1\ni s_0$ in $V$ such, that $f$
meromorphically extends onto $V_1\times W_1$.

\smallskip
\rm
\noindent
\bf Proof. \rm  Step 1. Let $s_0$ be a locally regular point of $S$. First of
all we remark that the family of analytic subsets $\{\Gamma _{f_s} \} _
{s\in S}$ is continuous at $s_0$ in $W_1\times X$. Indeed, let  $s_n\in S,
s_n\rightarrow
s_0 $ as $n\rightarrow \infty $. Then from (2.4.2) we see that $vol(\Gamma
_{f_{s_n}})$ are uniformly bounded and thus by {\sl Lemma 2.3.1} $\Gamma _
{f_{s_0}}
\subset (W_0 \times X)$ extends to a graph of meromorphic mapping over
$W$. This will be also denoted as $\Gamma _{f_{s_0}}$. Now
 if one could find a sequence
$s_n\in S, s_n \longrightarrow s_0$ as $n \longrightarrow \infty $ such
that $\Gamma _{f_{s_n}}\not\longrightarrow \Gamma _{f_{s_0}} $ in Hausdorff
 metric, then by Lemma 2.3.1, using the boundedness of volumes of
$\Gamma _{f_s}$,
one finds a subsequence, still denoted as $s_n$ such that $\Gamma _{f_{s_n}}
\longrightarrow \Gamma \supset  \Gamma _{f_{s_0}},$ but not equal
$\Gamma _{f_{s_0}}$.

But then, by the relations (2.3.2) and (2.3.3) of Lemma 2.3.1  one has
\smallskip
$$\lim_{n\longrightarrow \infty} {\sl vol}(\Gamma _{f_{s_n}})\ge \nu +
{\sl vol}(\Gamma _{f_{s_0}}), $$
\smallskip
\noindent
which contradicts (2.4.2) .

Let us prove now that the family ${\cal F} = \bigcup_{s\in S}\Gamma _{f_s}
\subset S\times W_0 \times X $ extends to a meromorphic family on $V_1
\times W_1\times X$ for some neighbourhood $V_1\ni s_0$.

Fix a point $z_0\in {\sl Reg}\Gamma _{f_{s_0}}\cap (W_0 \times X)$.
Here we
consider $\Gamma _{f_{s_0}}$ as analytic space itself. So ${\sl Reg}\Gamma_{f_{
s_0}}$ is connected dense subset in $\Gamma _{f_{s_0}}$ and $Sing\Gamma _
{f_{s_0}}:=\Gamma_{f_{s_0}}\setminus {\sl Reg}\Gamma_{f_{s_0}}$ is an analytic
subset of $\Gamma
_{f_{s_0}}$.

Now take a point $z_1\in {\sl Reg}\Gamma _{f_{s_0}}\cap (W\times X)$.
 Take a path $\gamma :
[0,1]
\longrightarrow {\sl Reg}\Gamma_{f_{s_0}}$ from $z_0$ to $z_1$. We shall prove
that there is a neighborhood $\Omega $ of $\gamma ([0,1])$ in $W\times X$
 and a neighborhood $V\ni s_0$ such that ${\cal F}\cap (V\times \Omega )$
 extends
to an analytic set in $V\times \Omega $.

By $T$ denote the set of those $t\in [0,1]$ that there exists a neighborhoods
 $\Omega_t \supset \gamma ([0,t])$ and $V_t\ni s_0$ such that ${\cal F}\cap V_t
 \times \Omega_t$
extends to an analytic set in $V_t\times \Omega_t$. Note that $T$ is open and
contains the origin. Now let $t_0$ be the cluster point of $T$. Find the chart
$\Sigma \cong \Delta ^q\times \Delta ^n$ for the space $W\times X$
in the neighborhood of $\gamma (t_0)$ with coordinates $u_1,\ldots,u_q,v_1,\ldots,v_n
$ in such a way that $\gamma (t_0)=0$ and $\Gamma _{f_{s_0}}\cap \Sigma =
\{ (u,v): v = F_0(u)\} $ for some holomorphic map $F_0:\Delta ^q\longrightarrow
 \Delta ^n$. By the Hausdorff continuity of our family $\{ \Gamma _{f_s}\}$ in
 $s_0$ $\Gamma _{f_s}\cap \Omega = \{ (u,v) : v = F_s(u)\} $ for $s$ close to
$s_0$, $F_s$ holomorphic and continuously depending on $s$.

Take $t_1\in T$ close to $t_0$, such that $\gamma ([0,t_1])\subset \Sigma $.
We have
some neighborhoods $V_{t_1}\ni s_0, \Omega_{t_1}\ni \gamma (t_1)$ such that ${\cal
F}$ extends analytically to $V_{t_1}\times \Omega_{t_1}$. Let $u_1\in \Delta^q  $ be
 such that $\gamma (t_1) = (u_1, F_0(u_1))$. Then there is a neighborhood, say
$\Delta ^q_r(u_1) $, such that $\Gamma _{f_s}\cap (V_{t_1}\times \Delta ^q_r
(u_1)
\times \Delta ^n)$ is defined by the  equation $v = F(u,s)$, where $F(u,s)=F_s
(u):V_{t_1}\times \Delta_r^q(u_1)\to \Delta^n$ as above. From the condition of
the {\sl Lemma}  we see that for $s\in S$ close to
$s_0$ $F(u,s)$ extends onto $\Delta^q$. So by {\sl Lemma 2.1.1} $F(u,s)$
extends
holomorphically to $V_{t_0,\varepsilon } \times \Delta^q_{1-\eps }$, where
$\varepsilon $ is arbitrarily small ($V_{t_0,\varepsilon }$ depending on
$\varepsilon $ ). But this means that ${\cal F}$ extends analytically onto
$V_{t_0}\times \Delta^q_{1-\varepsilon }\times \Delta ^n$. Thus $T$ is closed
and coinsides with $[0,1]$.

We proved in fact that for any  compact subset $R\subset {\sl Reg}\Gamma_
{f_{s_0}}
\cap (W_1\times X)$ there are neighborhoods $V_R\ni s_0,\Omega_R\supset R$
such that
$\Gamma _f$ analytically extend to $V_R\times \Omega_R$.

\noindent\sl
Step 2. \rm Cover the set ${\sl Sing}\Gamma _{f_{s_0}}\cap [\{ s_0\} \times
\bar W_1\times X
]$  with a finite number of open charts of the form  $ 1/2V_{\alpha }\times
\Omega_{\alpha }$, where $V_{\alpha }\cong \Delta ^q $ and $\Omega_{\alpha }\cong
\Delta ^n$, and such that $\Gamma _{f_{s_0}}\cap (V_{\alpha }\times
\Omega_{\alpha })$ is analytic cover of $V_{\alpha }$. By Step 1 we can find an open
neighborhoods $V_R\ni s_0$ and $\Omega_R\supset R = \Gamma _{f_{s_0}}\setminus
[\bigcup_{\alpha }1/2V_{\alpha }\times \Omega_{\alpha }]$ such that ${\cal F}$
analytically extends to $V_R\times \Omega_R$.

Fix now some $\alpha $. All that remained to prove is that $\Gamma_f$
analytically extends to $V'_{\alpha }\times V_{\alpha }\times \Omega_{\alpha }$ for
some neighborhood of $V'_{\alpha }\ni s_0$. But this again follows from
{\cal Lemma 2.1.1} applied to the coefficients of polynomials which define the
cover
$\Gamma_{f_s}\cap (V_{\alpha }\times \Omega_{\alpha })\to V_{\alpha }$.
\bigskip
\hfill{q.e.d.}
\medskip
 Divide the variables in $\cc^{n}=\cc^{p+q}$ into two groups
: $(z_1,...,z_p)$ and  $(u_1,...,u_q)$. Let $\Gamma $ an analytic set in
$\Delta^n\times X$. Fix some $0<r<1$ and put $\Gamma_z=\Gamma \cap (\{ z\}
\times \Delta^q\times
X)$. Consider a function
$$
v_{\Gamma }(z) = {\sl vol}_{2q}(\Gamma_z) = \int_{\Gamma_z}(dd^c\Vert u\Vert
^2 + w_h)^q.\eqno(2.4.3)
$$
\noindent
Here, as usually $w_h$ denotes the $(1,1)$-form canonically accosiated
to $h$.
This is well defined for $z\in \Delta^p$, for which $\dim (\Gamma_z)=q$.
Denote the set of such $z$ as $U$. $T:=\Delta^p\setminus U$ is contained
in at most countable union of locally closed proper analytic subvarieties
of $\Delta^p$, see [F]. We shall make use of the following result of
D. Barlet:

\it the function $v_{\Gamma }$ is locally bounded on $\Delta^p$,

\smallskip\noindent\rm
see [B] Th\'eor\`eme 3.
\smallskip\noindent\bf
Corollary 2.4.2. \it If $f:\Delta^n\to X$ is a meromorphic map then the Lelong
 numbers $\Theta((f^*w)^q,0)$ are finite for all $q=1,...,n$.
\smallskip\noindent\sl
Proof. \rm This is the same as to bound the Lelong numbers of the currents
$(dd^c
\Vert z\Vert^2+f^*w_h)^q$ Put $z'=(z_{i_1},...,z_{i_q})$ and $z^{''}=
(z_{j_1},...,z_{j_{n-q}})$. We have

$$
{1\over r^{2p}}\int_{\Delta_r^p\times \Delta_r^q}(f^*w+dd^c\Vert z'\Vert^2)^q
\wedge (dd^c\Vert z^{''}\Vert^2)^p={1\over r^{2p}}\int_{\Delta_r^p}(dd^c\Vert
z^{''}
\Vert^2)^p\int_{\Gamma_{z^{''}}}(dd^c\Vert z{'}\Vert^2+w)^q=
$$
$$
={1\over r^{2p}}\int_{\Delta_r^p}(dd^c\Vert z^{''}\Vert^2)^pv_{\Gamma }
(z^{''})\le
C\cdot { vol(\Delta_r^p)\over r^{2p}}\le C_1
$$
by the Barlet Theorem.
\smallskip
\hfill{q.e.d.}

\bigskip

\bigskip\noindent\sl
2.5. The Rothstein-type extension theorem and separate meromorphicity. \rm
\smallskip
Lemmas 2.3.1 and 2.4.1  together with
estimate in 3.2 will play a key role in the proof  of the  main result of
this paper. Hovewer we shall break here the main line of exposition to
show how they can be already applied to the classical question entitled
above.

We are going to treat the question of separate meromorhpicity in the form
proposed by Kazaryan and moreover, we are going to show the condition on the
image space $X$ to possess the mer.ext.property is "almost" not needed.
It was F. Hartogs who proved his
"Hartogs-type" extension theorem for holomorphic functions in order to
derive from this his separate holomophicity theorem. Rothstein did this
for meromorphic functions, starting from E. Levi extension theorem.
B. Shiffman
fulfilled
this program for meromorphic mappings. Namely he proved in [Sh-2] that if the
complex space $X$ possess a meromorphic extension property in dimension
$n$ then separately meromorphic mappings from $\Delta^n$ to $X$ are
meromorphic.

We shall prove here the Rothstein-type extension and separate meromorphicity
in the form of Siciak and Kazaryan. If the space $X$ posess a hol.ext.prop.
the separate holomorphicity in the form of Siciak was obtained, using the
approach of Shiffman, by O. Alehyane in [A-1].
\smallskip\noindent\bf
Corollary 2.5.1 \it Let $V\subset \cc^p$ and $W_0\subset\subset W\subset
\cc^q$
be a domains and let $E$ a nonpluripolar subset of $V$. Let further a
meromorphic mapping $f:V\times W_0\to X$ into a  complex space $X$ is given.
Suppose that for $z\in E$
the restriction $f_z:=f\mid_{\{ z\} \times W_0}$ is well defined and
extends meromrophically onto $\{ z\} \times W$. Then there is a pluripolar
subset $E'\subset E$ such that $f$ meromorphically extends onto a
neighborhood of $V\times W_0\cup (E\setminus E')\times W$.

\smallskip\noindent\sl
Proof. \rm
Take an exaustion $W_0\subset\subset W_1\subset\subset ... $ of $W$ by
relatively compact domains. For every $n$ we shall find a pluripolar subset
$E_n'\subset E$ such that $f$ extends meromorphically to the neighborhood
$P_n$ of $V\times W_0\cup (E\setminus E_n')\times W_n$. Then $P:=\bigcup_{
n\ge 1}W_n$ will be a neighborhood of $V\times W_0\cup (E\setminus
\bigcup_{n\ge 1}E_n')
\times W =V\times W_0\cup (E\setminus E')\times W$, where
$E':=\bigcup_{n\ge 1}E_n'$ is pluripolar in $V$ by the Josefson theorem.
Finally $f$ is extended to this neighborhood.

So, we can suppose additionaly that there is a domain $\tilde W\supset
\supset W$ such that for
any $z\in E$ $f_z$ is well defined and extends meromorphically onto
$cl(\tilde W)$-closure of $\tilde W$.

Next, fix some compact exhaustion $K_1\subset\subset ...\subset\subset K_n
\subset\subset ...$ of $X$. Denote now by $E_n$ the set of those $z\in E$
that $f_z(\tilde W)\subset K_n$. While $\bigcup_{(n)}E_n=E$, starting
from some $n_0$ all $E_n$ are not pluripolar. If we shall prove that there
exists a pluripolar $E_n'$ and extension of $f$ into the neighborhood
of $V\times W_0\cup (E_n\setminus E_n')\times W$, then again using Josefson
theorem we can finish the proof.

Thus we may additionally suppose that there is a compact $K\subset\subset X$
such that $f_z(\tilde W)\subset K$ for all $z\in E$.

Take $\nu $ as in Lemma 2.4.1 for $W_0\subset\subset W\subset\subset
\tilde W$ and $K$. Let $R$ be the maximal open
subset of $E$ such that $f$ extends to the neighborhood of $V\times W_0
\cup R\times W$. Put $S=E\setminus R$. Put further
$$
S_k=\{ z\in S: vol(\Gamma_{f_z})\le k{\nu \over 2}\} ,\eqno(2.5.1)
$$
\noindent
where graphs are taken over $\tilde W$. Note that by Lemma 2.3.1 $S_k$
are closed and they are increasing. Also $\bigcup_{k\ge 1}S_k=S$.
By Lemma 2.4.1
and maximality of $R=E\setminus S$ the sets $S_{k+1}\setminus S_k$ are
not locally regular at any of their points. Thus the sets
$S_1, S_2\setminus S_1,...$ are pluripolar and so is $S$.
\smallskip
\hfill{q.e.d.}
\smallskip\noindent\bf
Remarks. \rm 1. If $E=V$ and $X$ posseses a mer.ext. property
then we obtain the result of B. Shiffman, see [Sh-2].

\noindent 2. \rm If $X$ is compact and K\"ahler then it satisfies
this assumption, see [Iv-2].

\noindent 3. Let $X$ posseds a meromorphic extension property in dimension
$p+q$. Then $f$ extends onto the envelope of holomorphy of the neighborhood
of $V\times W_0\cup (E\setminus E')\times W$, which is obviously a
neighborhood of  $V\times W_0\cup E^*\times W$. Here we denote by $E^*$ the
set of locally regular points of $E$.

\noindent 4. Let us give the following
\smallskip\noindent\bf
Definition 2.5.1. \rm We shall say that meromorphic mappings into complex
space $X$ satisfy:

(h) a  Hartogs-type extension theorem in bidimension (p,q)
if every meromorphic map from
$$
H^p_q(r) = \{ (z_1,\ldots,z_{p+q})\in \cc^{p+q}:\Vert z'\Vert <r, \Vert z^{
''}\Vert
<1 \hbox{ or } \Vert z'\Vert <1, 1-r<\Vert z''\Vert <1\} , \eqno(2.5.2)
$$
\noindent
to $X$ extend meromorphically onto $\Delta^{p+q}$. Here $z'=(z_1,\ldots,z_p)$,
$z''=(z_{p+1},\ldots,z_{p+q})$;
\smallskip
(t) a strong Thullen-type extension theorem in bidimension (p,q) if
for any closed pluripolar subset $S$ of $\Delta^q$
every meromorphic map from

$$T^p_q(S,a,b):=(\Delta^q\times \Delta^p(a))\cup ((\Delta^q\setminus
S)\times \Delta^p(b)) \eqno(2.5.3)
$$
\noindent
to $X$ extends meromorphically onto $\Delta^{p+q}$, here $a<b$.
\smallskip\noindent
Note that Hartogs (p,q)-extendibility obviously implies the strong Thullen -
type one. Vice versa is not true. Example is given by a 3-fold constructed
by M. Kato in [Ka-2].

Namely Kato had constructed a compact three-fold $X$, which is a quotient of
$D=\{ [z_0:\ldots :z_3]\in \cc\pp^3:\vert z_0\vert^2+\vert z_1\vert^2<\vert z_2
\vert^2+\vert z_3\vert^2\} $ by a co-compact properly discontinuous subgroup
$G\subset {\sl Aut}(D)$. Denote by $\pi :D\to X$ the natural projection.
Consider a hyperplane $P=\{ z\in \cc\pp^3:z_0=0\} \cong \cc\pp^2$. Then
$D\cap P=\cc\pp^2\setminus \bar\bb^2$, here $\bb^2$ is a ball $\{ \vert z_1
\vert^2>\vert z_2\vert^2+\vert z_3\vert^2\} $ in $\cc\pp^3$. Thus
$\pi\mid_{D\cap P}:\cc\pp^2\setminus \bar\bb^2\to X$ cannot be extended onto
the neighborhood of any point on $\partial \bb^2$. So meromorphic mappings into
$X$ are not Hartogs (1,1)-extendable.

Let now $f:T^1_1(S,1/2,1)\to X$ be some meromorphic map. While
$f:T^1_1(S,1/2,1)$ is simply-connected we can consider a lifting $F=\pi^{-1}
\circ f: f:T^1_1(S,1/2,1)\to D\subset \cc\pp^3$. By a Thullen-type extension
theorem for meromorphic functions $F$ extends onto $\Delta^2$. But $F(\Delta^2)
\cap \partial D=\emptyset $ because one cannot touch $\partial D$ by bidisk.
Thus $\pi \circ F$ gives an extension of $f$ onto $\Delta^2$.

\smallskip
We have the following obvious corollary.
\smallskip\noindent\bf
Corollary 2.5.2. \it If in the conditions of the Corollary 2.5.1 a complex
space $X$
possess a strong Thullen-type extension property in bidimension $(p,q)$
and $E=V$ then $f$ extends meromorphically onto $V\times W$.

\smallskip\rm
Let us turn now to the separate meromorphicity.
\medskip\noindent\bf
Corollary 2.5.3. \it Let $E$ and $G$ be a nonpluripolar subsets in domains
$V\subset \cc^p$ and $W\subset \cc^q$ respectively. Let $F$ be some
pluripolar subset of in $V\times W$. Let further some mapping $f:E\times
G\setminus F\to X$ into a complex space $X$ is given. Suppose
that:

(i) for every $z\in E$, such that $\{ z\} \times G\not\subset F$ the
restriciton $f_z:=f\mid_{\{ z\} \times G}$ is well defined and
meromorphically extends meromorphically onto $\{ z\} \times G$, and

(ii) for every $w\in G$, such that $V\times \{ w\} \not\subset F$ the
restriciton $f^wz:=f\mid_{V\times \{ w\} }$ is well defined and
meromorphically extends meromorphically onto $V\times \{ w\} $.
\smallskip
Then there are pluripolar subsets $E'\subset E$, $G'\subset G$ and a
meromorphic mapping $\tilde f$ of some neighborhood of $(E\setminus E')
\times W\cup V\times (G\setminus G')$ into $X$ which extends $f$.
\smallskip\noindent\sl
Proof. \rm
Withought loss of generality as in the
proof of Rothstein-type theorem, we suppose that $f_z$ extends onto
$cl(\tilde W)$ for some $\tilde W\supset\supset W$ and for all $z\in E$,
and that there is a compact $K\subset\subset X$ such that $f_z(\tilde W)
\subset K$ for all $z\in E$.
The same for $f^w$-s.
\smallskip\noindent\sl
Step 1. \it There is a point $Z_0=(z_0,w_0)\in E\times G$ such that $f$
holomorphically extend to the neighborhood of $Z_0 $.
\smallskip\rm
Define
$$
E_k= \{ z\in E: vol(\Gamma_{f_z})\le k{\nu \over 2}\} . \eqno(2.5.4)
$$
\noindent
While $E$ is not pluriolar, there exists $k$ and $z_1\in E_{k+1}\setminus E_k$
such that $E_{k+1}\setminus E_k$ is locally regular at $z_1$. The
same reasoning as at the begining of the proof of Lemma 2.4.1 shows
that the family $\{ \Gamma_{f_z}:z\in E_{k+1}\setminus E_k\} $ is continuous
in the neighborhood of $z_1$.
Take a point $w_0\in W$ such that $f_{z_1}$ is \it holomorphic \rm
in the neighborhood of $w_0$. Remark that this $w_0$ can be taken to be
a locally regular point of $G_{l+1}\setminus G_l$ for some $l$.
From Hausdorff continuity of our family
in the the neighborhood of $z_1$ we get immediately that all $f_z$ are
holomorphic in the neighborhood of $w_0$ for $z\in E_{k+1}\setminus E_k$
close to $z_1$. Find a point $z_0$ close to $z_0$ where $f^{w_0}$ is
holomorphic. Now the separate analyticity theorem for functions tells us
thatthe point  $Z_0 = (z_0,w_0)$ is as needed.
\smallskip\noindent\sl
Step 2. \it End of the proof.
\smallskip\rm
Applying two times coordinatevise the Rothstein theorem we get the statement.
\smallskip
\hfill{q.e.d.}
\smallskip\noindent\bf
Remark. \rm

\noindent 1. Again, as above, if $X$ posseds a mer.ext.prop. then we can take
$E\setminus E'=E^*$ and $G\setminus G'=G^*$. This case was studied recently
in [A-2], using approach developped in [Sh-2] and [Sh-3].

\noindent 2. In the case $X=\cc\pp^1$ we obtain the theorem of Kazaryan, see
[Kz].

\bigskip\noindent
\bf 3 . Estimates of Lelong numbers from below.
\smallskip\noindent
\sl 3.1. Generalities on blowings-up.\rm
\medskip
First we recall the Hironaka Resolution Singularities Theorem. We shall
use the so called embedded resolution of singularities, see [H], [B-M].
Let us recall  the notion of the sequence of local blowings
up over a complex manifold $D$. Take a point $s_0\in D_0:=D$. Let $V_0$
be some neighbourhood of $s_0$ and $l_0$ smooth, closed submanifold
of $V_0$ of codimension at least two, passing through $s_0$.
Denote by $\pi _1 : D_1\longrightarrow
V_0$ the blowing up of $V_0$ along $l_0$. Call this \it a local blowing up
of $D_0$ along the locally closed center $l_0$. \rm  The exceptional divisor
$\pi ^{-1}(l_0)$ of this blowing up we denote by $E_1$.

We can repeat this procedure, taking a point $s_1\in D_1$, a neighborhood
$V_1$ of that point in $D_1$ and smooth closed submanifold  $l_1$ in $V_1$
of codimension at least two.
\smallskip\noindent\bf
Definition 3.1.1. \rm A finite sequence $\{ \pi ^j\} _{j=1}^N $ of such local
blowings up we call \it a sequence of local  blowings up over $s_0\in D$, \rm
 or a {\sl local regular modification}.

\smallskip
By $\{ l_j\} _{j=0}^{N-1}$ we denote
the corresponding centers in the neighborhoods $\{ V_j\}_{j=0}^{N-1} $ of
points
$\{ s_j\}_{j=0}^{N-1} $, and by  $\{ E_j\}_{j=1}^N $ the exceptional divisors,
$s_j\in D_j$.

If $V_j=D_j$ for all $j=0,...,N-1$, and points $s_j$ are not specified  we
call this sequence a sequence of (global) blowings up or \it a regular
modification. \rm In this case we put $\pi = \pi_1\circ ...\circ \pi_N$,
$\hat D=D_N$, $E$ denotes the exceptional divisor of $\pi $, i.e.
$E=\pi_N^{-1}
(l_{N-1}\cup ...\cup (\pi_1\circ ... \circ \pi_N)^{-1}(l_0)$.

\smallskip\noindent\bf
Theorem. \it Let $L$ be a subvariety of $D$. Then there exists a regular
(not local!) modification $\pi :\hat D\to D$ such that:

1) the strict transform $\hat L$ of $L$ is smooth;

2) $l_i\subset {\sl Sing}L_i$, where $L_i$ is a strict transform of $L$ by
$\pi_1\circ ...\circ \pi_j$.

\smallskip\noindent\rm See [H].

\smallskip
For the proof we refer to [B-M].

\smallskip\noindent\bf
Remark. \rm We shall use this Theorem only for the case {\sl dim }$D=3$ and
{\sl dim}$L=1$. Tn this case (while {\sl dim Sing}$L=0$), we need only
blowings up of points to resolve the singularities of $L$.
We shall need also the following three Lemmas about the behavior of meromorphic
mappings under the modifications. First let us introduce some more notations.
Let $f:D\setminus S\to X$ be a meromorphic map into a complex space $X$. Here
$D$ is a manifold and $S$ supposed to be closed and zero dimensional.

\smallskip\noindent\bf
Definition 3.1.2. \it Recall that a closed subset $S$ of a metric space is
called zero dimensional if
for any $s_0\in S$ and for almosr all $r>0$ the sphere centered at $s_0$
of radii $r$ do not intersect $S$.

\smallskip\rm
By $I(f)$
we shall denote the set of points of indeterminacy of $f$. By $I_p(f)$ those
components of $I(f)$ which have dimension at least $p$. By $I_{p,s}(f)$ the
set of components from $I_p(f)$ which pass through the point $s$; by $I_{p,R}
(f)$-the set of those components from $I_p(f)$ which intersect the set $R$.

Remark also that if $l$ is an (irreducible) analytic set in $D\setminus S$ of
pure dimension $p\ge 1$, then its closure is an (irreducible) analytic subset
of $D$, provided $S$ is zerodimensional.
\smallskip\noindent\bf
Lemma 3.1.1. \it  Let $f:B^n_*\to X$ be a holomorphic map of punctured
$n$-ball, $n\ge 2$, into a complex space $X$ which meromorphically extends
onto $B^n$. Suppose one can find a sequence $\{ \pi_j\}_{j=1}^{\infty }$
of blowings-up in such a way that:

(i) $\pi_1$ is a blow-up of $B_0^n:=B^n$ at $s_{0,1}=0$; $B^n_1:=\pi^{-1}_1
(B_0^n)$.

(ii) $\pi_{j+1}$ is a blow-up of $B_j^n$ at points $\{ s_{j,1},...,s_{j,N_j}\}
\subset \bigcup_{i=1}^{N_{j-1}}\pi^{-1}_j(s_{j-1,i})$; here $N_0=1$.

(iii) $f_j:=f\circ (\pi_1\circ ...\circ \pi_j)$ is holomorphic on
$B^n_j\setminus \{ s_{j,1},...,s_{j,N_j}\} $.

\noindent Then there exists $j_0$ such that $f_{j_0}$ is holomorphic on
$B^n_{j_0}$.
\smallskip\noindent\sl
Proof. \rm
By $\Gamma_f\subset B^n\times X$ denote the graph of $f$. Put $\Gamma = \Gamma
_f\cap (\{ 0\} )\times X)$. Write $\Gamma = \bigcup_{i=1}^N\Gamma_i$-
decomposition into irreducible components. Put $E_{j+1,i}=\pi_{j+1}^{-1}(s_{j,i
}), i=1,...,N_j, E_{j+1}=\bigcup_{i=1}^{N_j}E_{j+1,i}$ - the  exceptional
divisor of $\pi_{j+1}$. Usually we denote
by $E_{j+1,i}$ also all strict transforms of it by subsequent blowings-up.
\smallskip\noindent\sl
Step 1. For every $i=1,...,N$ there are a $j$ and $1\le k\le N_j$ such that
$f_{j+1}(E_{j+1,k})=\Gamma_i$.\rm
\smallskip
 We shall prove this for $i=1$. Take a point $a\in {\sl Reg}\Gamma_1\setminus
(\bigcup_{i=2}^N\Gamma_i)$. There is a holomorphic map $\phi :\Delta^2\to
\Gamma_f$, such that:

(i) $\phi (\Delta\times \{ 0\} )=a$,

(ii) $\phi (\Delta^2\setminus (\Delta\times \{ 0\} ))\subset \Gamma_f
\cap (B^n_*\times X)$,

(iii) for any $z'\in \Delta$ the map $ \phi_{z'}:=\phi \mid_{\{ z'\}
\times \Delta
}$  is proper and primitive (i.e. not multiple covered).

\smallskip
Such map can be constructed first as imbedding into the smooth manifold
$\tilde \Gamma_h$, which
is a modification of $\Gamma_f$. One should only chose $\phi $ to be
transversal
to the exceptional divisor of smoothing modification $\pi :\tilde \Gamma_f
\to \Gamma_f$
on each disk $\{ z'\} \times \Delta $, and then take $\pi \circ \phi $ as
$\phi $.

Put $\psi =p_1\circ \phi $, where $p_1:B^n\times X\to B^n$ is a natural
projection, and $\psi_{z'}:=\psi \mid_{\{ z'\} \times \Delta }$
for $z'\in \Delta $. By $\psi^j:\Delta^2\to B^n_j$  denote the lift
$(\pi_1\circ ...\circ \pi_j)^{-1}\circ \psi $ of $\psi $ by $\pi_1\circ ...
\circ \pi_j$.

There is a $j$ such that $\psi^j_{0'}$ is an imbedding of $\Delta $ (after
shrinking $\Delta $ if necessary). So (after shrinking of $\Delta$),
$\psi_{z'}$ is also imbedding for all $z'\in \Delta $. Remark that for all
$z'\in \Delta $ by construction one has $\lim_{z_k\rightarrow 0}(f_j\circ
\psi_{z'}^j)(z_k)=a$.
\smallskip
\noindent\sl Case 1. \rm $\psi^j_{0'}(0)\not\in \{ s_{j,1},...,s_{j,N_j}\} $
for
some $j$.

In this case clearly $f_j(E_{j,i})=\Gamma_1$, where $i\in \{ 1,...,N_j\} $ is
such that $E_{j,i}\ni \psi^j_{0'}(0)$.
\smallskip\noindent\sl Case 2. \rm $\psi^j_{0'}(0)=s_{j,i_j}$ for all $j$.

Then there exists $j_1$ such that $\psi^{j_1}_{z'}(0)\not=\psi^{j_1}_{0'}(0)$
 for $z'\not= 0$, but $\psi^{j_1-1}_{z'}(0)=\psi^{j_1-1}_{0'}(0)$. Now one
can take some $z'_0$ instead of $0'$ for which $\psi^{j_1}_{z'_0}\not\in
\{ s_{j_1,1},...,s_{j_1,N_{j_1}}\} $ and repeat the {\sl Case 1}.
\smallskip\noindent\sl
Step 2. \rm Let $\nu $ be from (2.4.1) for some fixed Hermitian metric $h$ on
$X$ and some compact $K\supset f(B^n_{1\over 2})$. Then by {\sl Lemma} 3.2.1
the sum of Lelong numbers of currents $f^*w+dd^c\Vert z\Vert^2,...,
(f^*w+dd^c\Vert z\Vert^2)^n$ is $\ge N\cdot \nu $.

If $f_{j_1}$ is not holomorphic say in $s_{j_1,1}$ we can take this point
instead of zero and repeat the {\sl Step 1}. If this procedure doesn't stop
then
the sum of  Lelong numbers must be infinite. This contradics the
meromorphicity of $f$ at zero, see {\sl Corollary 2.4.2}.
\smallskip
\hfill{q.e.d.}
\bigskip\rm
Let $f_0:D_0\setminus S_0 \to X$ be a meromorphic mapping into a
complex space $X$, $S_0$-being zerodimensional, $s_0\in S_0$, $\dim D_0=3$.
Put $ \{ l^{(0)}_1,...,l^{(0)}_{R_0} \} = I_{1,s_0}(f_0) $.

Suppose that $l^{(0)}_1$ is smooth. Consider the following sequence of
blowings-up.
\smallskip\noindent
1)  $p_1:D_1\to D_0$ is a blowing up of $D_0$ along $l^{(0)}_1$. By $E_1$
denote the exceptional divisor of $p_1$.

\rm  Suppose that $f_1:=f_0\circ p_1$ meromorphically extends onto $D_1
\setminus S_1$ with $S_1$ being zerodimensional. Take some $s_1\in p_1^{-1}
(s_0)\cap S_1$ (if such exists).  Let $p_1^{'}:D_1^{'}\to D_1$ be a regular
modification resolving the singularities of all noncompact components
of $I_{1,s_1}(f_1)$ which are contained in $E_1$. Denote them (and their
strict transforms) by $l^{(1)}_1,...,l^{(1)}_{P_1}$.
\smallskip\noindent
2)  $p_2:D_2\to D_1$ is the composition with $p_1^{'}$ of the regular
modification  $p_1^{''}:D_2\to D_1^{'}$ which is successive blown-up of
$l_1^{(1)}$, then $l_2^{(1)}$, and so on till $l_{P_1}^{(1)}$. By $E_2$
we denote the exceptional divisor of $p_1\circ p_2$.
\smallskip\noindent
n) Suppose that $p_n:D_n\to D_{n-1}$ is constructed.
\smallskip
By $E_n$ denote the exceptional divisor of $p_1\circ ...\circ p_n$.
Suppose that $f_n:=f_{n-1}\circ p_n$ meromorphically extends onto $D_n
\setminus S_n$ with $S_n$ being zerodimensional. Take some $s_n\in p_n^{-1}
(s_{n-1})\cap S_n$ (if such exists). Let $p_{n}^{'}:D_n^{'}\to D_n$ be a
regular all noncompact components
of $I_{1,s_n}(f_n)$ which are contained in $E_n$. Denote them (and their
strict transforms) by $l^{(n)}_1,...,l^{(n)}_{P_n}$.
\smallskip\noindent
n+1) $p_{n+1}:D_{n+1}\to D_n$ is the composition with $p_n^{'}$ of the
regular
modification  $p_n^{''}:D_{n+1}\to D_n^{'}$ which is successive blown-up of
$l_1^{(n)}$, then $l_2^{(n)}$, and so on till $l_{P_n}^{(n)}$. By $E_{n+1}$
we denote the exceptional divisor of $p_1\circ ...\circ p_{n+1}$.
\smallskip\noindent\bf
Lemma 3.1.2. \it  There exists a $n_0$ and $s_{n_0}\in S_{n_0}\cap
 (p_1\circ ...\circ p_{n_0})^{-1}(s_0)$ such that no noncompact  component of
$I_{1,s_{n_0}}(f_{n_0})$ is contained in $E_{n_0}$ - the exeptional divisor
of $p_1\circ ...\circ p_{n_0}$.
\smallskip\noindent\bf
Remark. \rm  This means that after applying $R_0$-times this {\sl Lemma} we
get a regular modification $p_{n_1}:(D_{n_1},S_{n_1},f_{n_1})\to (D_0,S_0,f_0)
$ and a point $s_{n_1}\in S_{n_1}\cap p_{n_1}^{-1}(s_0)$ such that all
components of $I_{1,s_{n_1}}(f_{n_1})$ are compact and contained in $p_{n_1}
^{-1}(s_0)$.
\smallskip\noindent\sl
Proof. \rm The proof follows from the previous {\sl Lemma 3.1.1} by sections.
 Take
 a point $a^{'}_0\in l^{(0)}_1\setminus S_0$. Find a coordinates $z_1,z_2,z_3$
in the polydisk neighborhood $\Delta^3$ of $a^{'}_0$ such that
\smallskip
(i) $a^{'}_0=0$,

(ii) $\Delta^3\cap S_0=\emptyset $,

(iii) $\Delta^3 \cap l^{(0)}_1=\{ z:z_1=z_2=0\} $.
\smallskip
For $z'\in \Delta^2_{z_2z_3}$ we put $\Delta_{z'}=\{ z'\} \times \Delta $
and  $f_{z'}:=f\mid_{\Delta_{z'}}$. Let $\nu_1$ be from {\sl
Lemma} 2.3.1. Put $\Sigma_0=\emptyset $ and
$$
\Sigma_j=\{ z'\in \Delta^2: {\sl vol}(\Gamma_{f_{z'}})\le {\nu_1\over 2
}\cdot j\} ,\hbox{ for } j\ge 1.\eqno(3.1.1)
$$
\noindent
From (2.3.2) we see that $\Sigma_j$ are closed, $\Sigma_j\subset \Sigma_{j+1}$,
 and $\bigcup_{j=1}^{\infty }\Sigma_j = \Delta^2$ . Find $j_1\ge 1$
such that $\Sigma_{j_1}\setminus \Sigma_{j_1-1}$ contains a
closed disk $\Delta(a^{''}_0)\subset l_1^{(0)}$ centered  at $a^{''}_0\in
l^{(0)}_1$.
Note that by {\sl Lemma} 2.3.1
 and the fact the $f$ \it is not holomorphic \rm in the neighborhood of $a^{''
}_0$ it follows that $\Sigma_{j_1}$  is contained in a proper analytic set
in the neighborhood of $a^{''}_0$ in $\Delta^2$. So we can find a disk
$\Delta (a_0)\subset  \Delta (a^{''}_0) $ such that for some $r_0>0$
$(\Delta_{r_0}\times
\Delta (a_0)\setminus  \{ 0\} \times \Delta (a_0))\subset \bigcup_{j>j_1}
\Sigma_j$. Remark that if $z_1^{'}
\in \{ 0\} \times \Delta(a_0)$ and $z_2^{'}\in (\Delta_{r_0}\times
\Delta (a_0)\setminus (\Delta(a_0)\times \{ 0\} )$ one has a jump of volume:
$\vert {\sl vol}(\Gamma_{f_{0,z_1^{'}}}) - {\sl vol}(\Gamma_{f_{0,z_2^{''}}})
\vert \ge \nu_1/2$.

 Suppose that for all $n$ the set $E_n$ contains a noncompact
component of $I(f_n)$, which projects by
$p_1\circ ...\circ p_n$ onto $l^{(0)}_1$ .
Repeating the arguments above we can construct a next sequence:

(i) $a_n$ a point on $l_1^{(n)}$;

(ii) a disk $\Delta(a_n)$ and disk $\Delta_{r_n}$ such that if $z_1^{'}
\in \{ 0\} \times \Delta(a_n)$ and $z_2^{'}\in  \Delta_{r_n})\times
\Delta(a_n))\setminus (\{ 0\} \times \Delta(a_n))$ one has a jump of volume:
$\vert {\sl vol}(\Gamma_{f_{n,z_1^{'}}}) - {\sl vol}(\Gamma_{f_{n,z_2^{''}}})
\vert \ge \nu_1/2$.

(iii) $p_n(\Delta(a_n)\subset\subset \Delta(a_{n-1})$.
\smallskip
\noindent
Let $a\in \bigcap_{n=1}^{\infty }(p_1\circ ... \circ p_n)(\Delta(a_n))$.
Take a two-ball $B_a^2$ centered at $a$ and transversal to $l_1^{(0)}$.
 Note that if this ball was chosen sufficiently small that $f\mid_{B^2
_a}$ is meromorphic and by properties (i)-(iii) above all $f\circ (p_n\circ
...\circ p_1)$ are \it essentially meromorphic \rm i.e. not holomorphic. This
contradics {\sl Lemma} 3.1.1.
\smallskip
\hfill{q.e.d.}
\smallskip
 In the sequel we shall repeatedly use the following statement.
\smallskip\noindent\bf
Lemma 3.1.3. \it Let $f:H^{n+1}(r)\to X$ be a meromorphic map into a complex
space $X$, which possess a meromorphic extension property
in dimension $n$. Then $f$ meromorphically extends onto $\Delta^{n+1}\setminus
 S $, where $S=S_1\times ...\times S_{n+1}$ and all $S_j\subset \Delta $ are
of harmonic measure zero.
\smallskip\noindent\sl
Proof. \rm We shall prove that $f$ extends meromorphically onto $E^{n+1}(r)
\setminus S$ with $S$ as described. Here
$$
E^{n+1}(r) = (\Delta^{n-1}_r\times \Delta_r\times \Delta )\cup (\Delta_r^{n-1}
\times \Delta \times A_{1-r,1}) = \bigcup_{z_1\in \Delta_r} E^n_{z_1}(r)
$$
\noindent
in notations of 2.2. This will clearly imply the statement of {\sl Lemma}.

According to the assumption of $n$-dimensional extension property of $X$,
for every $z_1\in \Delta_r$ the restriction $f_{z_1}=f\mid_{E^n_{z_1}}(r)$
extends meromorphically onto $D_{z_1}:=\{ z_1\} \times \Delta_r^{n-2}
\times \Delta^2$.

After shrinking, we can suppose that all $f_{z_1}$ are meromorphic in the
neighborhood of $\bar D_{z_1}$ in $\{ z_1\} \times \cc^n$. Put $D:=\Delta_r
^{n-2}\times \Delta^2$. Take $\eps >0$ and consider $D_{\eps }= \{ z\in D:
d(z,\partial D)\ge \eps \} $. Denote by $\Omega_{\eps }$ the maximal open
subset of $\Delta_r$ such that $f$ meromorphically extends onto $\Omega
\times D_{\eps }$. Let $S^{\eps }=\Delta_r\setminus \Omega_{\eps }$. Let
$\nu_p$ be as in (5) of {\sl Lemma 2.1.1}. Put $\nu ={\sl inf} \{ \nu_p
\cdot \eps^{2(n-p)}:p=0,...,n-1\} $.

Consider a following closed subsets of $S^{\eps }$: $S^{\eps }_j=\{ z_1\in S^{
\eps }: {\sl vol}(\Gamma_{f_{z_1}})\le {\nu \over 2}\cdot j\} $. Note that
$S^{\eps }_{j+1}\supset S^{\eps }_j$ and $S^{\eps }=\cup_{(j)}S^{\eps }_j$.
$S^{\eps }_{j+1}\setminus S^{\eps }_j$ is of harmonic measure zero in
$\Delta_r\setminus S^{\eps }_j$. Really, would $s_0\in S^{\eps }_{j+1}
\setminus S^{\eps }_j$ be some regular point of $S^{\eps }_{j+1}$, then by
{\sl Lemma 2.4.1} $f$ would meromorphically extend onto $V\times D$ fore
some neighborhood $V\ni s_0$. This contradics the maximality of
$\Omega_{\eps }$. By the Josefson $S^{\eps }$ is polar.

Further $S_1=\cup_{\eps }S^{\eps }$ is polar, so $f$ extends onto
$(\Delta_r\setminus S_1)\times D$. Repeating the same arguments for other
coordinates we obtain the statement of {\sl Lemma}.
\smallskip
\hfill{q.e.d.}

Consider a
meromorphic map $f_0:D_0\setminus S_0\to X$ into a complex space which
possesses a meromorphic extension property in dimension $n$, $S_0$ is
zerodimensional and ${\sl dim }D_0=n+1$. Let $\pi_1:D_1\to D_0$ be some
regular modification.
\smallskip\noindent\bf
Lemma 3.1.4. \it $f_0\circ \pi_1$ extends to a meromorphic map $f_1:D_1
\setminus S_1\to X$, where $S_1$ is zero dimensional, closed and pluripolar
subset of $D_1$.
\smallskip\rm
The proof, which is similar to that of previous {\sl Lemma}, will be omited.

We shall also make use of the fact that the set of harmonic measure zero
on the plain has zero dimension, see [Gl].
\smallskip\noindent
\rm Note that {\sl Lemma 3.1.3} gives part (i) of {\sl Theorem 2}.

\bigskip\noindent
\sl 3.2. Estimates.
\smallskip\rm
Let $X$ be a complex space, equipped with some Hermitian metric $h$. By
$\omega $ we denote, as usually, (1,1)-form canonically associated with
$h$.
Let $S_0$ be some zero dimensional closed subset of a complex manifold  $D_0$
and let $f:D_0\setminus S_0 \longrightarrow X$ be some meromorphic map such
that $cl(f(D_0\setminus S_0))\subset K $-some compact in $X$. Suppose that some
sequence $\{ \pi_j\} _{j=1}^N$ of local blowings-up $\pi_j:D_j\to D_{j-1}$
over the point $s_0\in S_0$
is given. Denote by $f_j$
the lifting of $f_{j-1}$ onto $D_j$, i.e. $f_j=f_{j-1}\circ \pi _j$, where
$f_0=f$. Suppose now that all $f_j$ extend meromorphically onto $D_j\setminus
S_j$, where $S_j$ are closed zero dimensional subsets of $D_j$. As before, by
$l_{j-1}$ we denote the center of $ \pi_j $ and by $E_j=\pi_j^{-1}(l_{j-1})$
the corresponding exeptional divisor (and all its strict transforms under
$\pi_{j+1},...,\pi_N$). Denote by $\hat S = S_N$ and by $\hat D = D_N$.
Put also $\pi =\pi_1\circ ...\circ ^pi_N:\hat D\to D_0$.
We suppose moreover that $f_j|_{E_j\setminus S_j}$
meromorphically extends onto $E_j$. For every $j=1,...,N$ by $r_j$ denote the
rank of $f_j|_{E_j}$.

\smallskip\noindent\bf
Remark. \rm  In what follows the locall blowings up will
appear in the following context. Take a point $s_j\in E_j$ and let $l_j$ be
the smooth \it compact \rm component of $I_{1,s_j}(f_j)$ (note that ${\sl
codim }
l_j\ge 2$). In this
 case we shall take $V_j=D_j$ and blow up $D_j$ along $l_j$ to obtain  an
exeptional divisor $E_{j+1}$. If $l_j$ was compact then $E_{j+1}$ is also
compact. Such configurations will contribute to our estimates of Lelong
numbers.

Consider now currents $T_p = (f^*\omega + dd^c\Vert z\Vert ^2)^p$ on
$D_0\setminus (S_0\cup I(f_0))$. More accurately the currents $T_p$
one can define as $T_p=((p_1\mid_{\Gamma_f})_*p_2^*w + dd^c\Vert z\Vert
^2)^p$.
For each $T_p$ we consider the Lelong
number of
$T_p$ in $s_0$:
\smallskip
$$
\Theta (T_p,s_0) = \lim_{\varepsilon \searrow 0}\sup 1/\varepsilon
^{2(n-p)}\int_{B_{s_0}(\varepsilon )\setminus (S_0\cup I(f_0))} T_p \wedge
(dd^c\Vert z \Vert ^2)^{2(n-p)},
\eqno(3.2.1)
$$
\smallskip\noindent
which can of course take infinite value, $n={\sl dim}D_0 $. By $\nu_p=\nu_p(K)$
  we
denote,
as before the minima of volumes of $p$-dimensional compact subvarieties
in $X$,
which are contained in $K$. By $\sigma_p$ denote the number of those
$j$ that $l_j$ is compact and $r_{j+1}=p$.
\smallskip
\noindent
\bf Lemma 3.2.1. \it There is a constant $C=C(h,K)$, depending only on
Hertmitian
 metric $h$ and compact $K$, such that for any meromorphic map $f$ as above
one has the next estimate for the Lelong numbers of the currents $T_p$ in
the point $s_0\in S$:
\smallskip
$$
\Theta (T_p,s_0) \ge C\cdot \sigma _p\cdot \nu_p \eqno(3.2.2)
$$
\smallskip
\noindent
\rm Proof. Take only thouse $E_{j_1}$ which are wlowings-up of
{\it compact} $l_{j_1-1}$ and ${\sl rk}f_{j_1}\mid_{E_{j_1}}=p$. To simplify
the notations in what follows we drop the subindice $1$ in $j$-s.

Take an $\eps $-neighborhood $V^{\eps }$ of $\bigcup_{j=1}^{\sigma_p}E_j$ with
respect to some metrik on $\hat D$. Remove from $V^{\eps }$ the $\eps $-
neighborhood of intersections $E_i\cap E_j, i\not= j$,  $\eps $-neighborhood
 of $\hat S$ and $\eps $-neighborhood of $I(f_N)$. For $\eps >0$ small enough
we obtain the union $\bigsqcup _{j=1}
^{\sigma_p}V_j^{\eps }$ of pairwise disjoint open sets: $\bigsqcup_{j=1}
^{\sigma_p}V_j^{\eps } = V^{\eps }\setminus ((\bigcup_{i\not= j}E_i\cap E_j)^
{\eps }\cup \hat S^{\eps }\cup I(f_N)^{\eps }) $. By $\bar V_j^{\eps }$
denote
the closure
 of $V_j^{\eps }$ in $\hat D$.

Denote by $W_j^{\eps }$ the image of $\bar V_j^{\eps }$ under the blown-down
mapping $\pi :\hat D\to D_0$. Starting from here we put $s_0=0$.
Note that for $\eps >0$
suffuciently small $W_j^{\eps }\cap
W_i^{\eps } = \{ 0\} $ for $i\not= j$. Note further that $\pi^{-1}\mid_{W_j
^{\eps }\setminus \{ 0\} } : W_j^{\eps }\setminus \{ 0\} \longrightarrow
\hat D$ is correctly defined and, in fact, $\pi^{-1}(W_j^{\eps }\setminus \{ 0
\} ) = \bar V_j^{\eps }\setminus E_j$. Denote by $F_j^{\eps }$ the closure
in $D_0\times X$ of the graph of $f\mid_{W_j^{\eps }\setminus \{ 0\} } $.

 So $F^{\eps }
_j\cap (\{ 0\} \times X)$ is compact in $\{ 0\} \times X$ of finite
Hausdorff $p$-measure. In fact $F^{\eps }_j\cap (\{ 0\} \times X)$ is the
closure of the image
of $f_N(E_j\setminus (S_j^{\eps }\cup \bigcup_{i\not= j}E^{\eps }_i)\cup I(f_N)
^{\eps })$, so $F^{\eps }
_j\cap (\{ 0\} \times X) $ is contained in a compact subvariety $A_j=f(E_j)$
 of $\{ 0\} \times X$ of complex dimension $p$.

Let $K$ be a compact in $X$ containing $cl(f(D_0\setminus S_0))$. Cover $K$ by
a finite number of open sets $\{ U_
{\alpha }\} $, which are biholomorphic to an analytic subsets of $\Delta^{m_1}
_{\alpha }$-
unit polydisk in $\cc^{m_1}, m_1\ge m={\sl dim}X  $. Appropriate coordinate
functions on $U_{\alpha }$
are denoted by
$u_1^{\alpha },...,u_{m_1}^{\alpha }$. In each $\Delta_{\alpha }^{m_1}$ find
an
orthonormal
basis $e_{\alpha }^1,...,e_{\alpha }^{m_1}$ of $(1,0)$-forms with respect to
the
Hermitian metric $h_{\alpha }$ on $\Delta_{\alpha }^{m_1}$, so that the form
$w_h$ can be written as $w_h =
\sum_{j=1}^{m_1}{i\over 2}e_{\alpha }^j\wedge e_{\alpha }^{\bar j}$. Now put
$\Vert u^{\alpha }\Vert^2=\sum_{j=1}^{m_1}\vert u_j^{\alpha}\vert^2$ and
express
$$
{i\over 2}\partial \bar\partial \Vert u^{\alpha }\Vert^2 = \sum_{j,k=1}^{m_1}
c^{\alpha }_{j\bar k}{i\over 2}e_{\alpha }^j\wedge \bar e_{\alpha }^k, \eqno
(3.2.3)
$$
\noindent
where $[c^{\alpha }_{j\bar k}]_{j,k=1}^{m_1}$ is strictly positive-definite
matrix depending on $u_{\alpha }\in U_{\alpha }$. Find constants $C^1_{\alpha }
,C^2_{\alpha }$ such that
$$
C^1_{\alpha }\cdot w_h \le {i\over 2}\partial \bar\partial \Vert u_{\alpha }
\Vert^2=\sum_{j,k=1}^{m_1}c^{\alpha }_{j\bar k}{i\over 2}e^j_{\alpha }\wedge
\bar e^k_{\alpha }\le
$$

$$
\le C^2_{\alpha }\cdot \sum_{j=1}^{m_1}{i\over 2}e^j_{\alpha }\wedge \bar
e^k_{\alpha } = C^2_{\alpha }w_h. \eqno(3.2.4)
$$
\noindent
Denote by
$A_j^{-\eps }:={\sl Reg}A_j\setminus (f_j(E_j\cap S_j^{\eps }\cap \bigcup_
{i\not= j}
E_i^{\eps })))$. Cover $A_j^{-\eps }$ by a finite number of open sets $V_i$
satisfying the following conditions:
\smallskip\it
(i)  $V_i\subset \subset U_{\alpha }$ for some $\alpha $. \rm

If $\phi_{\alpha }:U_{\alpha }\to \Delta_{\alpha }^{m_1}$ is coordinate
imbedding, and $\phi_{\alpha }\mid_{V_i}$ is a proper imbedding into some
open subset $V^{'}_i$ of $\Delta_{\alpha }^{m_1}$, then on $V^{'}_i$ one can
introduce coordinates $u^i_1,...,u^i_{m_1}$ such that \it

(ii) $V^{'}_i$ is polydisk in coordinates $u^i$ and $\phi_{\alpha }(V_i\cap
A_j^{-\eps }) = \{ u^i: u^i_{p+1}=...=u^i_{m_1}=0\}  $.

(iii) ${1\over 2}\cdot dd^c\Vert u^{\alpha }\Vert^2\le dd^c\Vert u^i\Vert^2\le
 2\cdot dd^c\Vert u^{\alpha }\Vert^2$.

(i{\sl v}) Each point of $A_j^{-\eps }$ belongs not more than to $2p+1$ of
$V_i$.
\smallskip
\noindent\rm
If $p=m_1$, then the second condition is not needed. Put $U=V^{'}_i$ to
simplify the notations and denote by $B^n_r$ the ball of radiu
$r$ centered at zero, $0<r<\eps $, where $\eps $ is
choosen suffuciently small to garantee that the restriction onto $F_j^{\eps }
\cap (B^n_r\times U)$ of the projection ${\sl pr}:B^n\times U
\longrightarrow B^n\times U_{u_1,...,u_p}$ is proper. Here $U=U_{u_1,...,
u_p}\times U_{u_{p+1},...,u_{m_1}}$. We had droped the indice $i$ in $u^i$.

 Now put $C(K)=\max_{\alpha }\{  C^1_{\alpha },C^2_{\alpha },2 \}  $. Further
put $Z_j =
{\sl pr}(F_j^{\eps }\cap (B^n\times U)$. It is a closed subvariety of
dimension $n$ in $B^n\times U_{u_1,...,u_p}$.

Consider the following Hermitian metric on $B^n\times U_{u_1,...,u_p}$:
$e = \sum_{j=1}^n{i\over 2}dz_j\otimes d\bar z_j + \sum_{k=1}^p{i\over 2}du_k
\otimes d\bar u_k $. Denote by ${\sl vol}_eZ_j$ the volume of $Z_j$ with
respect to $e$. Then putting $\Vert u\Vert_1^2 = \sum_{j=1}^p\vert u_j\vert^2$
, we have
$$
{\sl vol}_e(Z_j) \le \int_{(B^n_r\times U)\cap F_j^{\eps }}(dd^c(\Vert z
\Vert^2 + \Vert u\Vert_1^2))^n\le
$$
$$\le {2^n\over n!}\int_{(B^n_r\times U)\cap F_j^{\eps }}(dd^c(\Vert z
\Vert^2 + \Vert u\Vert_1^2))^p\wedge
(dd^c(\Vert z\Vert)^2)^{n-p}\le
$$
$$\le C^p(K) {2^n\over n!}\int_{(B^n_r\times U)\cap F_j^{\eps }}
(dd^c\Vert z\Vert^2 + w)^p\wedge
(dd^c\Vert z\Vert^2)^{n-p} =
$$
$$ = C(K)^p{2^n\over n!}\int_{B^n_r\cap W_j^{\eps }}T_p\wedge (dd^c\Vert z
\Vert^2)^{n-p} .
\eqno(3.2.5)
$$
\smallskip
 From
the well known lower bound of volumes of analytic varietes, see [A-T-U], one
has
$$
{\sl vol}_{2n}(Z_j)\ge C\cdot r^{2n-2p}\cdot {\sl vol}_{2p}(Z_j^0).\eqno(3.2.
6)
$$
\noindent here $Z_j^0 = Z_j\cap (\{ 0\} \times U) = A_j^{-\eps }\cap U$ and
$C$ is absolute constant. Now
from
(3.2.5) and (3.2.6) we get that
$$
C^p(K){2^n\over n!}\int_{B^n_r\times \cap W_j^{\eps }}T_p\wedge (d\bar d^c
(\Vert z\Vert )^2)^{n-p}\ge {\sl vol}_e(Z_j) \ge
$$
$$
\ge C\cdot r^{2n-2p}{\sl vol}_e(Z_j^0) = r^{2n-2p}C\cdot {\sl vol}_e
(A_j\cap U) \ge r^{2n-2p}C\cdot 1/C^p(K)\cdot {\sl vol}_h(A_j\cap U).
$$
\noindent So

$$
{1\over r^{2n-2p}}\int_{B^n_r\cap W_j^{\eps }}T_p\wedge (dd^c\Vert z\Vert^
2)
^{2n-2p}\ge
 {n!C^{2p}(K)\over 2^n}\sum_i{\sl vol}_h(A^{-\eps }_j\cap V_i) \ge
$$

$$
\ge {n!C^{2p}(K)\over 2^n(2p+1)}{\sl vol}_h(A^{-\eps }_j).
$$

\noindent
Remained to take sum over $j=1,...,\sigma_p$ and let $\eps \to 0$.

\smallskip

\hfill{q.e.d.}

\bigskip\noindent\sl
Proof of Theorem 1.\rm
\smallskip\noindent
(a) is proved in {\sl Lemma 3.1.3}.

\noindent
(b) Take a ball $D_0$ centered at $s_0$ and such that all
components of $I_1(f_0)$ intersecting $D_0$ pass through $s_0$. Denote
$S_0=D_0\cap S$ and put
$f_0=f\mid_{D_0\setminus S_0}$.  Denote by $l_{1}^{(0)},...,l_{N_0}^{(0)}$
all
components of codimension two in $I_{s_0}(f_0)$. In the sequel we shall say
that a regular modification $\pi :(D_1,S_1,f_1)\to (D_0,S_0,f_0)$ is given,
having in mind that $f_0$ is defined and meromorphic on $D_0\setminus S_0$,
where $S_0$ is zero dimensional, and that $f_1=f_0\circ \pi $ extends
meromorphically onto $D_1\setminus S_1$, where $S_1$ again is zero
dimensional by Lemma 3.1.4. We denote by $I_R(f)$ the set of all irreducible
components of $I(f)$ which intersect the subset $R$.

\smallskip\noindent\sl
Step 1. \it There is a regular modification $\pi_1:(D_1,S_1,f_1)\to (D_0,S_0,
f_0)$, such that the set  $I_{\pi_1^{-1}(s_0)}(f_1)$  doesn't contain
components of codimension
two, which $\pi_1$ maps onto some of  $l_{1}^{(0)},...,l_{N_0}^{(0)}$ .
\smallskip\noindent\sl
Proof of Step 1. \rm Let $p_1:(D_1,S_1,f_1)\to (D_0,S_0,f_0)$ be a regular
modification which is a composition of resolution of singularities of $l_1^
{(0)}$ and blow-up of strict transform of $l_1^{(0)}$. By $l^{(1)}=\{
l_1^{(1)},...,l_{N_1}^{(1)}\} $ denote the set of all components of $I_{p_1^
{-1}(s_0)}(f_1)$ whom $p_1$ projects onto $l_1^{(0)}$.

Further, let $p_2:(D_2,S_2,f_2)\to (D_1,S_1,f_1)$ be a regular modification,
which is a composition of resolution of singularities of all $l_1^{(1)},...,
l_{N_1}^{(1)}$ and successive blowings-up of their strict transforms, and
$f_2=f_1\circ p_2$.

Suppose $p_n:(D_n,S_n,f_n)\to (D_{n-1},S_{n-1},f_{n-1})$ is constructed.
Let us denote by  $l^{(n)}=\{ l_1^{(n)},...,l_{N_n}^{(n)}\} $ the set of all
components of $I_{1,(p_1\circ ...\circ p_n)^{-1}}(f_n)$ whom $p_n\circ ...
\circ p_1$ projects onto
$l_1^{(0)}$. Let $p_{n+1}:
(D_{n+1},S_{n+1},f_{n+1})\to (D_n,S_n,f_n)$ be a regular modification,
which is a composition of resolution of singularities of all $l_1^{(n)},...,
l_{N_n}^{(n)}$ and successive blowing-up of their strict transforms. Again
$f_{n+1}:=f_n\circ p_{n+1}$.

{\sl Lemma} 3.1.2 tells us that for some $n_0$ we have $l^{(n_0)}=\emptyset $.
 Denote by $l_1^{(n_0)},...,l_{N_{n_0}}^{(n_0)}$ the set of all
components of $I_{(p_1\circ ...\circ p_{n_0})^{-1}}(f_{n_0}) $ whom
$p_1\circ ...\circ p_{n_0}$ projects onto some of
$l_1^{(0)},...,l_{N_0}^{(0)}$. We proved that $(p_1\circ ...
\circ p_{n_0})(\bigcup_{i=1}^{N_{n_0}}l_i^{(n_0)})\subset \bigcup_{i=2}^{N_0}
l_i^{(0)}$.

By $E_1$ denote a component of the exeptional divisor of $p_1\circ ...\circ
 p_{n_0}$ such that $(p_1\circ ...\circ  p_{n_0})(E_1)=l_1^{(0)}$.

Repeating this procedure for the strict transform  of $l_2^{(0)}$ instead of
$l_1^{(0)}$ and so on till $l_{N_0}^{(0)}$ we get the regular modification
 $\pi_1:(D_1,S_1,f_1)\to (D_0,S_0,f_0)$ such that no $1$-dimensional
component of $I(f_1)$ is mapped by $\pi_1 $ onto some of
$l_{1}^{(0)},...,l_{N_0}^{(0)}$. By $E_j$ denote the union of the components
of the exeptional
 divisor of $\pi_1$ such that $\pi_1(E_j)=l_j^{(0)}$. {\sl Lemma} 3.1.2 gives
us now the statement of Step 1.

By  $l^{(1)}=l_1^{(1)},...,l_{N_1}^{(1)}$ denote the set of all one-dimensional
 componets of $I_{\pi_1^{-1}(s_0)}(f_1)$. Note that they all are contained in
the fiber $\pi_1^{-1}(s_0)$. In particular they are compact.
\smallskip\noindent\sl
Step 2. \it There is a regular modification $\pi_2:(D_2,S_2,f_2)\to (D_1,S_1,
f_1)$, such that the set $I_{(\pi_1\circ \pi_2)^{-1}(s_0)}(f_2)$  doesn't
contain any component of codimension two. Moreover

\noindent $\Theta_{tl}(f^*w,s_0)\ge N_1\cdot \nu $. Here the total Lelong
number is defined as $\Theta_{tl}(f^*w,s_0) :=$

\noindent $ = \Theta (f^*w,s_0) + \Theta ((f^*w)^2,s_0)$ and $N_1$ is a
number of  1-dimensional components of

\noindent $I_{\pi_1^{-1}(s_0)}(f)$.
\smallskip\noindent\sl
Proof of Step 2. \rm Construct inductively the following sequence of regular
modifications:
\smallskip\noindent
(1) $p_2:(D_2,S_2,f_2)\to (D_1,S_1,f_1)$ is a composition of resolution of
singularities of all $l^{(1)}_{1},...,l_{N_1}^{(1)}$
 and successive blowings-up of their strict transforms.
\smallskip\noindent
Suppose $p_n:(D_n,S_n,f_n)\to (D_{n-1},S_{n-1},f_{n-1})$ is constructed.
Let us denote by  $l^{(n)}=\{ l_1^{(n)},...,l_{N_n}^{(n)}\} $ the set of all
one-dimensional components of $I_{(\pi_1\circ p_2\circ ...\circ p_n)^{-1}}
(f_n)$
whom
$\pi_1\circ p_2\circ ...\circ p_n$ projects onto some of
 $l_1^{(0)},...,l_{N_1}^{(0)}$.
\smallskip\noindent
(n+1) $p_{n+1}:(D_{n+1},S_{n+1},f_{n+1})\to (D_n,S_n,f_n)$ is a regular
modification, which is a composition of resolution of singularities of all
$l_1^{(n)},...,l_{N_n}^{(n)}$ and successive blowing-up of their strict
transforms.
\smallskip
By {\sl Lemma }3.1.2 for some $n_1$ we have $l^{(n_1)}=\emptyset $. This
means that there are compact divisors $E_1^{(1)},...,E_{N_1}^{(1)}$ such that
$p^{n_1}:=\pi_1\circ p_2 \circ ...\circ p_{n_1}$ maps $E_i^{(1)}$ onto
$l_i^{(1)
}$ and $f_{n_1}\mid_{E_i^{(1)}\setminus S_{n_1}}$ meromorphically extends onto
$E_i^{(1)}$. {\sl Lemma }3.2.1 gives us $\Theta_{tl}(f^*w,s_0)\ge N_1\cdot
\nu $. Put $\pi_2 := p_2\circ ...\circ p_{n_1}, D_2:=D_{n_1}, f_2:=f_{n_1},
 S_2:=S_{n_1}$.

If $I_{(\pi_1\circ \pi_2)^{-1}(s_0)}(f_{n_1})\not= \emptyset $ we can repeat
this procedure
for $f_{n_1}$ instead of $f_1$ and obtain the following sequence of  regular
modifications:

$$
...\buildrel \pi_{n+1} \over\to D_n \to ...\to D_2 \buildrel \pi_2 \over\to
D_1 \buildrel \pi_1 \over\to D_0
$$
Here, if we put $\pi^n=\pi_1\circ ...\circ \pi_n$ and
 $l^{(n)}=\{ l_1^{(n)},...,l_{N_{n}}^{(n)}\} =I_{1,(\pi^n)^{-1}(s_0)}(f_n)$,
then $I_{1,(\pi^{n+1})^{-1}(s_0)}(f_{n+1})$ doesn't contain components which
$\pi^{n+1}$ maps onto some of $l_1^{(n)},...$ $...,l_{N_{n}}^{(n)}$. So if
$I_{1,(\pi^n)^{-1}(s_0)}(f_n)$ is not empty for all $n$ then
 $\Theta_{tl}(f^*w,s_0)\ge \Sigma_{i=1}^nN_i\cdot \nu $ for all $n$. This
contradics the condition of the Theorem. Thus $I_{(\pi_1\circ \pi_2)^{-1}(s_0
)}(f_{n_1})=\emptyset $ for some $n$, and consequently $l^{(n)}=\emptyset $
for some $n$. Hence the {\sl Step} 2 is proved.
\smallskip\noindent\sl
Step 3. \it There is a regular modification  $\pi_3:(D_3,S_3,f_3)\to (D_0,S_0,
f_0)$, such that  $f_3$ is holomorphic in the neighborhood of $\pi_3^{-1}
(s_0)$.

\noindent\sl
Proof of Step 3. \rm Let  $\pi_1\circ \pi_2:(D_2,S_2,f_2)\to (D_0,S_0,f_0)$
be a regular
modification constructed in {\sl Step} 2. If $S_2\cap I_{(\pi_1\circ \pi_2)^
{-1}(s_0)}(f_2)$ is not empty then take a point $s_2\in
S_2\cap I_{(\pi_1\circ \pi_2)^{-1}(s_0)}(f_2)$.
If $\pi_1\circ \pi_2$ is not biholomorphic,
then there is a compact curve $C\ni s_2$, $C\subset (\pi_1\circ \pi_2)^{-1}
(s_0)$, $C\cong
\cc\pp^1$. If $\pi_1\circ \pi_2$ is biholomorphic take $C=\{ z_2=z_3=0\} $ -
line in
$D_0=B^3$ with $s_0\in C$.

By $B^2$ denote some two-ball transversal to $C$ at $s_2$.

\smallskip\noindent\sl
Case 1.  $f_2\mid_{B^2
\setminus S_2}$ holomorphically extends onto $B^2$.

\rm In this case  $f_2$
is holomorphic in the neighborhood of $s_2$. To see this, take a Stein
neighborhood $W$ of $\Gamma_{f(\bar B^2)}$ in $D_2\times X$. Then the needed
statement follows from the Hartogs extension theorem for holomorphic functions.
But this contradicts the choice of $s_2$.

\smallskip\noindent\sl
Case 2. $f_2\mid_{B^2}$ is not holomorphic (i.e. essentially meromorphic).

\rm Denote
by $p_3:\hat B^2\to B^2$ the finite sequence of blowings-up of points such
that $f_2\circ p_3:\hat B^2\to X$ is holomorphic. Take a regular modification
$\pi_3:D_3\to D_2$ as an extension of $p_3$ along $C$. Let $p_3'$
(correspondingly $\pi_3^{'}$ be the last blowing-up in the sequence $p_3$
(corr.$\pi_3$). Let $l_3$ (corr.$E_3$) denote the exeptional divisor of $p^{'}
_3$ (corr.$\pi_3^{'}$). Then either $(f_2\circ \pi_3)\mid_{E_3}$ is nonconstant
 or $l_3\in I_{1,\pi_3^{-1}(s_0)}(f_2\circ \pi_3)$. In both cases we have a
contribution to $\Theta_{tl}(f^*w,s_0)$ of at least $\nu $. If $f_2\circ \pi_3
$ is not holomorphic in the neighborhood of $(\pi_2\circ \pi_3)^{-1}(s_0)$ then
 take $\pi_4$ to be a composition of a regular modification described in
{\sl Step 2} and just above. We get a sequence of regular modifications

$$
...\buildrel \pi_{n+1} \over\to D_n\to ...\to D_2 \buildrel \pi_2 \over\to D_1
\buildrel \pi_1 \over\to D_0
$$
such that $\Theta_{tl}(f^*w,s_0)\ge n\cdot \nu $. So it is finite by the
assumption of the Theorem. That means that for $f_n\circ \pi_n$ with $n$
big enough only the first case can occur. Thus $f_n\circ \pi_n$  is
holomorphic in the neighborhood of $(\pi_1\circ ...
\circ \pi_n)^{-1}(s_0)$.

The proof of {\sl Step 3} and thus of {\sl Theorem 1} is completed.
\smallskip
\hfill{q.e.d.}

\spaceskip=4pt plus3.5pt minus 1.5pt
\spaceskip=5pt plus4pt minus 2pt
\font\csc=cmcsc10
\font\tenmsb=msbm10
\def\rr{\hbox{\tenmsb R}}
\def\cc{\hbox{\tenmsb C}}
\newdimen\length
\newdimen\lleftskip
\lleftskip=2.5\parindent
\length=\hsize \advance\length-\lleftskip
\def\entry#1#2#3#4\par{\parshape=2  0pt  \hsize%
\lleftskip \length%
\noindent\hbox to \lleftskip%
{\bf[#1]\hfill}{\csc{#2 }}{\sl{#3}}#4%
\medskip
}
\ifx \twelvebf\undefined \font\twelvebf=cmbx12\fi

\bigskip\bigskip
\bigskip\bigskip
\centerline{\twelvebf References.}
\bigskip


\entry{A-1}{Alehyane O.:}{Une extension du th\'eor\`eme de Hartogs pour les
applications s\'epar\'ement holomorphes.} C.R.Acad.Sci. Paris, t.323,
S\'erie I (1996).

\entry{A-2}{Alehyane O.:}{Applications s\'epar\'ement m\'eromorphes dans
les espaces analytique.} Preprint.

\entry{A-T-U}{Alexander H., Taylor B., Ullman J.:}{Areas of Projections of
Analytic Sets.} Invent.\ math., {\bf16}, 335-341, (1972).

\entry{Ba}{Barlet D.:}{Majoration du volume des finres g\'en\'eriques et
la forme g\'eom\`etrique du th\'eor\`eme d'aplatissement.}
 Semiar P. Lelong - H. Skoda, Lect. Notes Math., {\bf 822}, 1-17, (1978/79).

\entry{B-T}{Bedford E., Taylor B.:}{ A new capacity for plurisubharmonic
functions.}  Acta\ Math. (1982) {\bf149}, 1-40, (1982).

\entry{B-M}{Bierstone, E., Milman, P.D.:}{Local resolution of singularities.}
 Proc.\ Symp.\ Pure.\ Math. {\bf52}, 42-64, (1991).


\entry{Gl}{Golysin G.:}{Geometric function theory.} Nauka,\ Moscow (1966).

\entry{H-S}{Harvey R., Shiffman B.:}{A characterization of holomorphic chains
.} Ann.\ Math. {\bf99}, 553-587, (1974).

\entry{Hi}{Hironaka H.:}{Introduction to the theory of infinitely near
singular points.} Mem.\ Mat.\ Inst.\ Jorge \ Juan, {\bf28}, Consejo Superior
de Investigaciones Cientificas, Madrid (1974).

\entry{Hs}{Hirschowitz A.:}{Les deux types de m\'eromorphie diff\'erent.}
J.\ reine \& \ ang. \ Math., {\bf313} , 157 - 160, (1980).

\entry{F}{Frisch J.}{Points de platitude d'une morphisme d'espaces
analytiques complexes.} Invent. math. {\bf4} , 118 - 138, (1967).

\entry{Iv-1}{Ivashkovich S.:}{The Hartogs phenomenon for holomorphically
convex K\"ahler manifolds} Math. USSR Izvestiya {\bf29} N 1, 225-232 (1987).

\entry{Iv-2}{Ivashkovich S.:}{The Hartogs-type extension theorem for the
meromorphic maps into compact K\"ahler manifolds.} Invent.\ math. {\bf109}
 , 47-54, (1992).

\entry{Ka-1}{Kato, M.:}{Compact complex manifolds containing "global
spherical shells".} Proc.\ Int.\ Symp.\ Alg.\ Geom.\ Kyoto. 45-84, (1977).

\entry{Ka-2}{Kato, M.:}{Examples on an Extension Problem of Holomorphic Maps
and Holomorphic 1-Dimensional Foliations.} Tokyo\ Journal\ Math. {\bf13}, n 1,
 139-146, (1990).

\entry{Kz}{Kazaryan M.:}{Meromorphic extension with respect to
grouops of variables.} Math. USSR Sbornik {\bf53}, n 2, 385-398, (1986).

\entry{Kl}{Klimek M.:}{Pluripozential theory.} London.\ Math.\ Soc.\ Monographs
 ,\ New\ Series 6, (1991).


\entry{Os}{Osgood W.:}{Lehrbuch der Funktionstheorie, Bd.II.1.} Leipzig.
 \ Teubner (1929).

\entry{Re}{Remmert R.:}{Holomorphe und meromorphe Abbildungen komplexer
R\"aume.} Math.\ Ann. {\bf133}, 328-370, (1957).


\entry{Sa}{Sadullaev A.:}{Plurisubharmonic measures and capacities on
complex manifolds.} Russian\ Math.\ Surv. {\bf36}, No 4 , 61-119, (1981).

\entry{Sh-1}{Shiffman B.:}{Extension of Holomorphic Maps into Hermitian
Manifolds.} Math.\ Ann. {\bf194}, 249-258, (1971).

\entry{Sh-2}{Shiffman B.:}{Separately meromorphic mappings into compact
K\"ahler manifolds.} Contributions to Complex Analysis and Analytic Geometry
 (H. Skoda and J.-M.Tr\'epreau, eds) Vieweg, Braunschweig, Germany, 243-250
(1994).

\entry{Sh-3}{Shiffman B.:}{Hartogs theorems for separately holomorphic
mappings into complex spaces.} C.R.Acad.Sci. Paris, t.310, S\'erie I, 89-84,
 (1990).

\entry{Si-1}{Siu Y.-T.:}{Every stein subvariety admits a stein neighborhood.}
 Invent.\ Math. {\bf38}, N 1, 89-100, (1976).



\entry{St}{Stolzenberg G.:}{Volumes, Limits and Extension of Analytic
Varieties.} Springer \ Verlag (1966).

\bigskip\bigskip
Universit\'e de Lille-I

UFR Math\'ematiques

59655 Villeneuve d'Ascq Cedex

France
\smallskip
ivachkov@gat.univ-lille1.fr

\bigskip\bigskip
IAPMM Acad Sci. of Ukraine

Naukova 3/b, 290053 Lviv

Ukraine
\end

\bye